\documentclass[11pt]{article}

\newcommand{\tn}{\textnormal}
\newcommand{\mc}{\mathcal}
\newcommand{\mb}{\mathbb} \newcommand{\lp}{\left(}
\newcommand{\rp}{\right)} \newcommand{\lb}{\left\lbrace}
\newcommand{\rb}{\right\rbrace} \newcommand{\la}{\left\langle}
\newcommand{\ra}{\right\rangle}
\newcommand{\ls}{\left[}
\newcommand{\rs}{\right]}	
\usepackage[normalem]{ulem}
\usepackage {epsfig,amsmath,euscript}

\usepackage[latin1]{inputenc}

\usepackage[english]{babel}
\usepackage{color}
\usepackage{subcaption}
\usepackage{amsmath}
\usepackage{amsthm}

\usepackage{amssymb}

\usepackage{ae}
\usepackage{float}
\usepackage[T1]{fontenc}
\usepackage{graphicx}
\usepackage{epstopdf}
\usepackage{longtable}
\usepackage[dvipsnames]{xcolor}
\definecolor{rred}{rgb}{0.7,0.0,0.2}
\definecolor{bblue}{rgb}{0.2,0.0,0.7}

\newcommand{\secref}[1]{Section \ref{sec:#1}}

\newtheorem{definition}{Definition}
\newtheorem{proposition}{Proposition}
\newtheorem{lemma}{Lemma}
\newtheorem{remark}{Remark}

\newcommand{\seclab}[1]{\label{sec:#1}}

\newcommand{\eqlab}[1]{\label{eq:#1}}
\renewcommand{\eqref}[1]{(\ref{eq:#1})}

\newcommand{\figref}[1]{Fig.~\ref{fig:#1}}
\newcommand{\figlab}[1]{\label{fig:#1}}

\newcommand{\corref}[1]{Corollary~\ref{cor:#1}}
\newcommand{\corlab}[1]{\label{cor:#1}}
\newcommand{\defnref}[1]{Definition~\ref{definition:#1}}
\newcommand{\defnlab}[1]{\label{definition:#1}}
\newcommand{\lemmaref}[1]{Lemma~\ref{lemma:#1}}
\newcommand{\lemmalab}[1]{\label{lemma:#1}}
\newcommand{\remref}[1]{Remark~\ref{remark:#1}}
\newcommand{\remlab}[1]{\label{remark:#1}}

\newtheorem{cor}{Corollary}

\makeatletter
\newcommand{\tpitchfork}{%
	\vbox{
		\baselineskip\z@skip
		\lineskip-.52ex
		\lineskiplimit\maxdimen
		\m@th
		\ialign{##\crcr\hidewidth\smash{$-$}\hidewidth\crcr$\pitchfork$\crcr}
	}%
}
\makeatother
\usepackage{tikz}

\usetikzlibrary{positioning}
\usetikzlibrary{arrows}
\usetikzlibrary{arrows.meta}
\usetikzlibrary{calc}

\title{Complex oscillatory dynamics in a three-timescale El Ni\~no Southern Oscillation model}

\author{Panagiotis Kaklamanos$ ^1 $ and Nikola Popovi\'c$ ^1 $}
\date{$ ^1 $Maxwell Institute for Mathematical Sciences and School of Mathematics,\\ University of Edinburgh}

\begin{document}

\maketitle

\pagestyle{myheadings}
\thispagestyle{plain}

\begin{abstract}
We study the three-timescale dynamics of a model that describes the El Ni\~no Southern Oscillation (ENSO) phenomenon, which was proposed in [A. Roberts, J. Guckenheimer, E. Widiasih, A. Timmermann, and C. K. Jones,
Mixed-mode oscillations of El Ni\~no--Southern Oscillation, Journal of the Atmospheric Sciences,
73 (2016), pp. 1755--1766]. While ENSO phenomena are inherently characterised by the presence of multiple distinct timescales, the above model has previously been studied in a two-timescale context only. Here, we uncover the geometric mechanisms that are responsible for complex oscillatory dynamics in a three-timescale regime, and we demonstrate that the system exhibits a variety of qualitatively different behaviours, such as mixed-mode oscillation (MMO) with ``plateaus" -- trajectories where epochs of quiescence alternate with dramatic excursions -- and relaxation oscillation. The latter, although emergent also in the two-timescale context in appropriate parameter regimes, had not been documented previously. Moreover, we show that these mechanisms are relevant to models from other fields of ecological and population dynamics, as the underlying geometry is similar to the unfolding of Rosenzweig--MacArthur-type models in three dimensions.
\end{abstract}

\section{Introduction}
The El Ni\~no Southern Oscillation (ENSO) phenomenon is associated with a variation in winds and sea surface temperatures over the Pacific Ocean, due to large-scale interactions between the ocean and the overlying
atmosphere \cite{dijkstra2005nonlinear}. It is composed of two phases: the El Ni\~no (``little boy'') phase, when warm water and weak trade winds develop in the East-Central equatorial Pacific Ocean, and the La Ni\~na (``little girl'') phase, when cold water and strong trade winds occur in the East-Central equatorial Pacific Ocean \cite{dijkstra2005nonlinear}. Although highly irregular, these patterns are oscillatory in nature; correspondingly, various multiscale systems of differential equations have been proposed to model the ENSO phenomenon: thus, for instance, it has been argued that the oscillations observed therein are a product of Hopf bifurcations, and that irregularities arise through Shilnikov-type homoclinic orbits and homoclinic--heteroclinic dynamics in the governing ordinary differential equations (ODEs) \cite{dijkstra2005nonlinear,jin1998simple,roberts2016mixed,timmermann2003nonlinear}.

In the present paper, we consider, in particular, the following model from \cite{jin1998simple,roberts2016mixed, timmermann2003nonlinear}:
\begin{subequations}\eqlab{rob}
\begin{align}
\frac{\tn{d}T_1}{\tn{d}t}&=-\alpha\lp T_1-T_r\rp-\varepsilon \mu \lp T_2-T_1\rp^2, \eqlab{rob-a}\\
\frac{\tn{d}T_2}{\tn{d}t}&=-\alpha\lp T_2-T_r\rp+\zeta \mu \lp T_2-T_1\rp\ls T_2-T_{\rm sub}\lp T_1, T_2, h_1\rp\rs, \eqlab{rob-b}\\
\frac{\tn{d}h_1}{\tn{d}t}&=r\ls-h_1-\frac{bL\mu(T_2-T_1)}{2\beta}\rs, \eqlab{rob-c}
\end{align}
\end{subequations}
where 
\begin{equation}
T_{\rm sub}(T_1, T_2, h_1)=\frac{T_r+T_0}{2}-\frac{T_r-T_{r_0} }{2}\tanh{\ls\frac{H-z_0+h_1+bL\mu(T_2-T_1)/\beta}{h^*}\rs}. \eqlab{Tsub} 
\end{equation}
Here, the variable $T_1$ corresponds to the equatorial temperature of the Western Pacific Ocean; $T_2$ corresponds to the equatorial temperature of the Eastern Pacific Ocean; and the variable $h_1$ denotes the thermocline depth of the Western Pacific. The first terms in the $(T_1, T_2)$-subsystem, Equations~\eqref{rob-a} and \eqref{rob-b},  represent the tendency of the system towards a climatological mean state $T_r$ in the absence of ocean dynamics; that is, without any interaction terms, the temperatures $T_1$ and $T_2$ would converge to the mean state $T_r$. The nonlinear interaction terms for the ocean dynamics in these $(T_1,T_2)$-equations depend on the temperature difference between $T_1$ and $T_2$, as well as on the difference of $T_2$ and the subsurface temperature $T_{\rm sub}$ given by \eqref{Tsub} -- the parameter $T_{r_0}$ therein corresponds to a mean Eastern equatorial
temperature, attained at a depth of about $75$ metres. Finally, the $h_1$-equation, Equation~\eqref{rob-c}, describes the tendency of the thermocline depth towards mean climatological conditions; correspondingly, the damping parameter $r$ is associated with the characteristic timescale of this process. 

Introducing the change of variables 
\begin{equation}
S=T_2-T_1, \quad T=T_1-T_r, \quad\text{and}\quad h=h_1+k,
\end{equation}
as well as the transformation
\begin{equation}
x=\frac{S}{S_0},\quad y=\frac{T}{T_0},\quad z=\frac{h}{h_0},\quad\text{and}\quad \tau_1=\frac{t}{t_0},
\end{equation}
where $S_0$, $T_0$, $h_0$, and $t_0$ are suitably chosen reference values, together with an appropriate scaling of the remaining parameters in \eqref{rob},  we can non-dimensionalise the governing equations as follows:
\begin{subequations}\eqlab{elno}
	\begin{align}
	x'&=x\big[ x+y+c(1-\tanh{(x+z)})\big]+\rho\delta(x^2-ax), \eqlab{elno-a}\\
	y'&=-\rho\delta(ay+x^2),\eqlab{elno-b}\\
	z'&=\delta(k-z-\tfrac{x}{2}); \eqlab{elno-c}
	\end{align}
\end{subequations}
see \cite{roberts2016mixed} for details.
For the variables and parameters in \eqref{elno}, we now have 
\begin{gather}
 \begin{gathered}
 x\leq0, \quad  y\in \mb{R},\quad z\geq0,\\ c\in(1,c_0),\quad  k\in (0,1),\quad  a\in (0,a_0),\quad\text{and}\quad 0<\delta,\rho\ll1,
 \end{gathered}\eqlab{consts}
\end{gather}
for some fixed $c_0>1$ and $a_0>0$.
We note that the first two equations in \eqref{elno} have been somewhat rearranged compared to \eqref{rob}, and we reiterate that the variable $ x $ corresponds to the temperature difference between the Eastern and Western Pacific surface water; $ y $ corresponds to the departure of the Western Pacific surface ocean temperature from some reference mean temperature; and $ z $ represents the Western Pacific thermocline depth anomaly. Therefore, it is apparent that the parameters $c$, $k$, $a$, $\delta$, and $\rho$ are associated to the rates of the aforementioned processes, and that they can hence be traced back to Equation~\eqref{rob}. 
\begin{remark}
    We restrict our analysis to the scenario where $c_0$ and $a_0$ in \eqref{consts} are not so large as to render the $c(1-\tanh{(x+z)})$-term dominant in Equation~\eqref{elno-a}, nor such that $\rho\delta a = \mc{O}(1)$ in Equations~\eqref{elno-a} and \eqref{elno-b}. 
    \remlab{caun}
\end{remark}

In particular, the bursting behaviour of ENSO can be associated with mixed-mode oscillatory (MMO) dynamics in \eqref{elno}, i.e. with patterns that feature alternating epochs of perturbed slow dynamics, such as epochs of quiescence or small-amplitude oscillation (SAO), followed by large, relaxation-type excursions or large-amplitude oscillation (LAO); cf.~\figref{outline} for examples. Such dynamics appear frequently in singularly perturbed systems with multiple timescales; while they are fairly well understood in systems with two timescales \cite{brons2006mixed,desroches2012mixed,wechselberger2005existence}, progress has been made only over the past few years on mixed-mode dynamics in a three-timescale context \cite{jalics2010mixed,kaklamanos2022bifurcations,krupa2008mixed,krupa2012mixed, letson2017analysis}.

It is widely accepted, from a climatological viewpoint, that processes in environmental and meteorological phenomena occur on various timescales \cite{dijkstra2005nonlinear}; in particular, in \cite{roberts2016mixed}, the model in Equation~\eqref{elno} was noted to exhibit dynamics on three distinct timescales when both $\delta$ and $\rho$ are positive and small. In this work, we therefore extend the analysis of \eqref{elno} from \cite{roberts2016mixed}, which was restricted to a two-timescale context, i.e. to $ \delta>0 $ sufficiently small and $ \rho = \mc{O}(1) $ in \eqref{elno}, to the three-timescale regime where both $ \delta $ and $ \rho>0 $ are small. We apply Geometric Singular Perturbation Theory (GSPT) \cite{cardin2017fenichel,fenichel1979geometric,kaklamanos2022bifurcations,letson2017analysis} to analyse the behaviour of \eqref{elno} upon variation of the (positive) parameters $c$, $k$, and $a$. Specifically, we focus on the various types of complex oscillatory dynamics exhibited by the system, and on the geometric mechanisms that underlie transitions between qualitatively different scenarios; cf. again \figref{outline}. 

\begin{figure}[ht!]
    \begin{subfigure}[b]{0.48\textwidth}
	\centering
	\begin{tikzpicture}
	\node at (0,0){
	\includegraphics[scale=0.25]{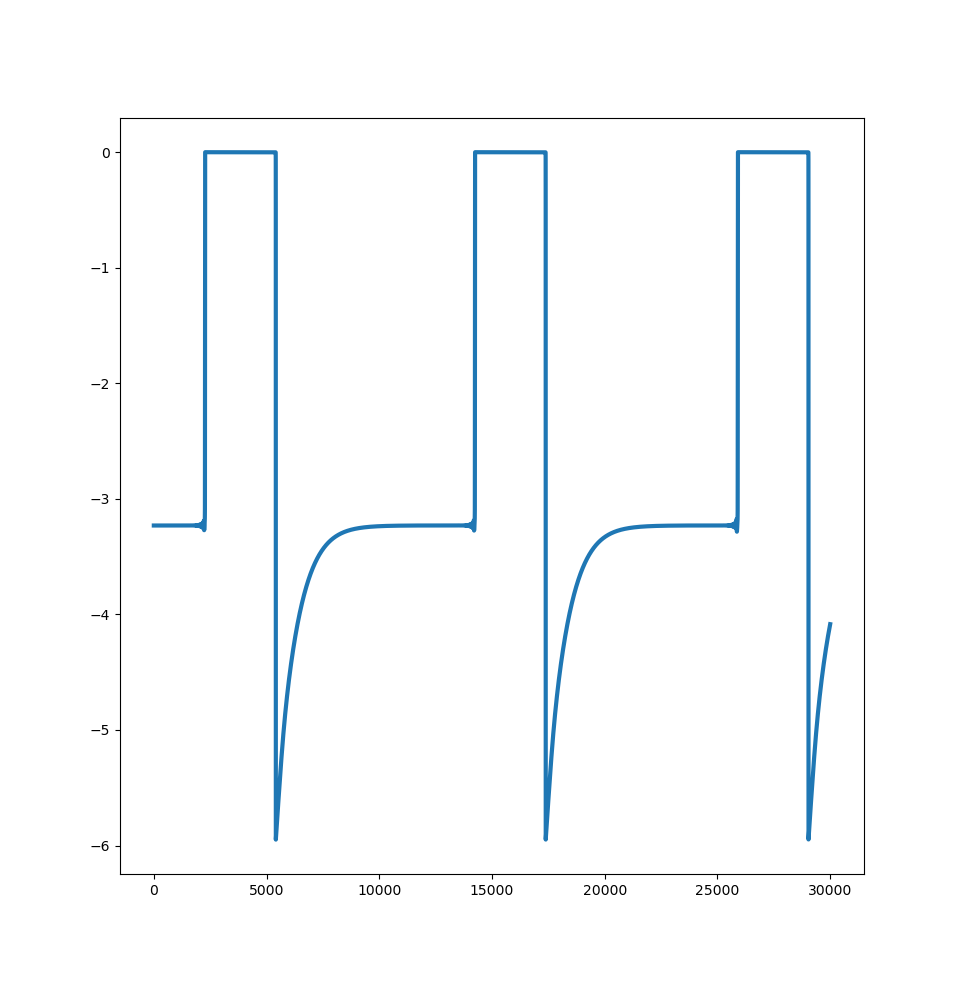}};
	
	\node at (0.15,-2.9) {$t$};
	\node at (-2.9,0) {$x$};
	
	\end{tikzpicture}
		\caption{$c = 3.75$, $k=0.34$, $a=2.8$}
	\end{subfigure}
	\begin{subfigure}[b]{0.48\textwidth}
	\centering
	\begin{tikzpicture}
	\node at (0,0){
	\includegraphics[scale=0.25]{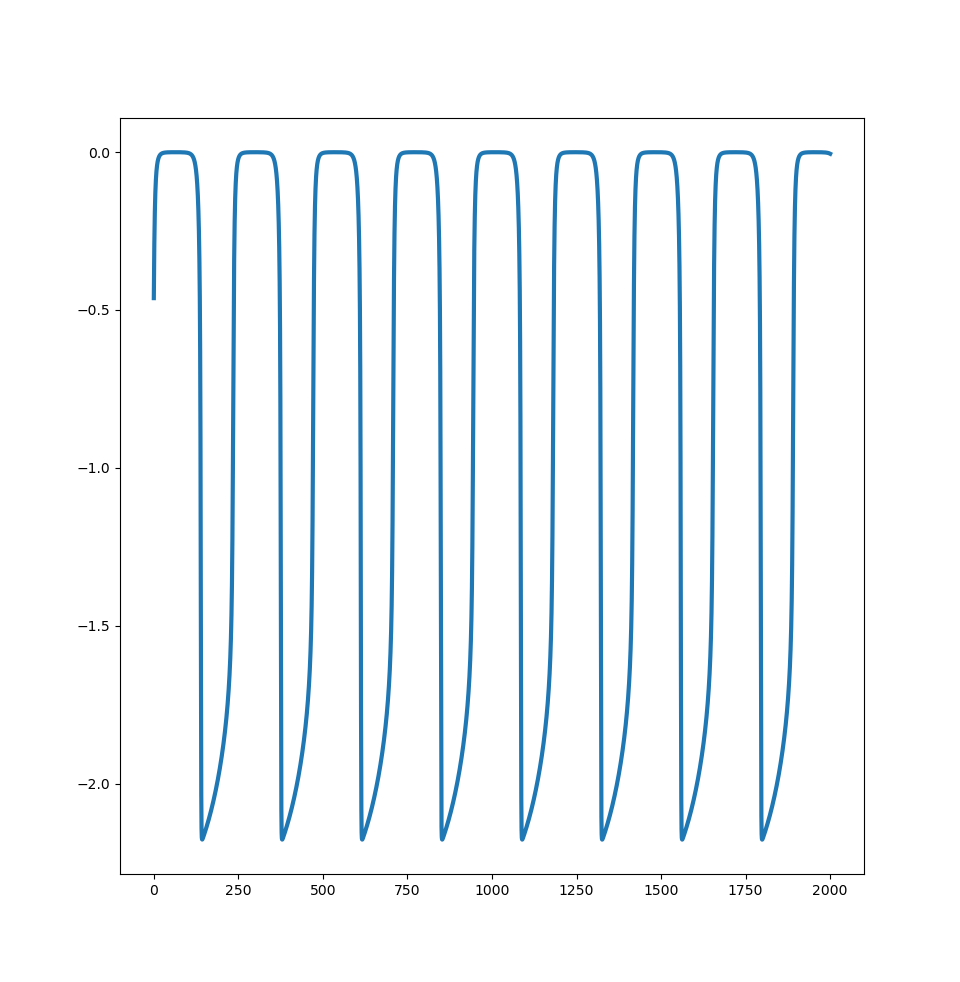}};
	
	\node at (0.15,-2.9) {$t$};
	\node at (-2.9,0) {$x$};
	
	\end{tikzpicture}
	\caption{$c = 1.5$, $k=0.34$, $a=0.5$} 
	\end{subfigure}
		\\
	\centering
	\begin{subfigure}[b]{0.48\textwidth}
	\centering
	\begin{tikzpicture}
	\node at (0,0){
	\includegraphics[scale=0.25]{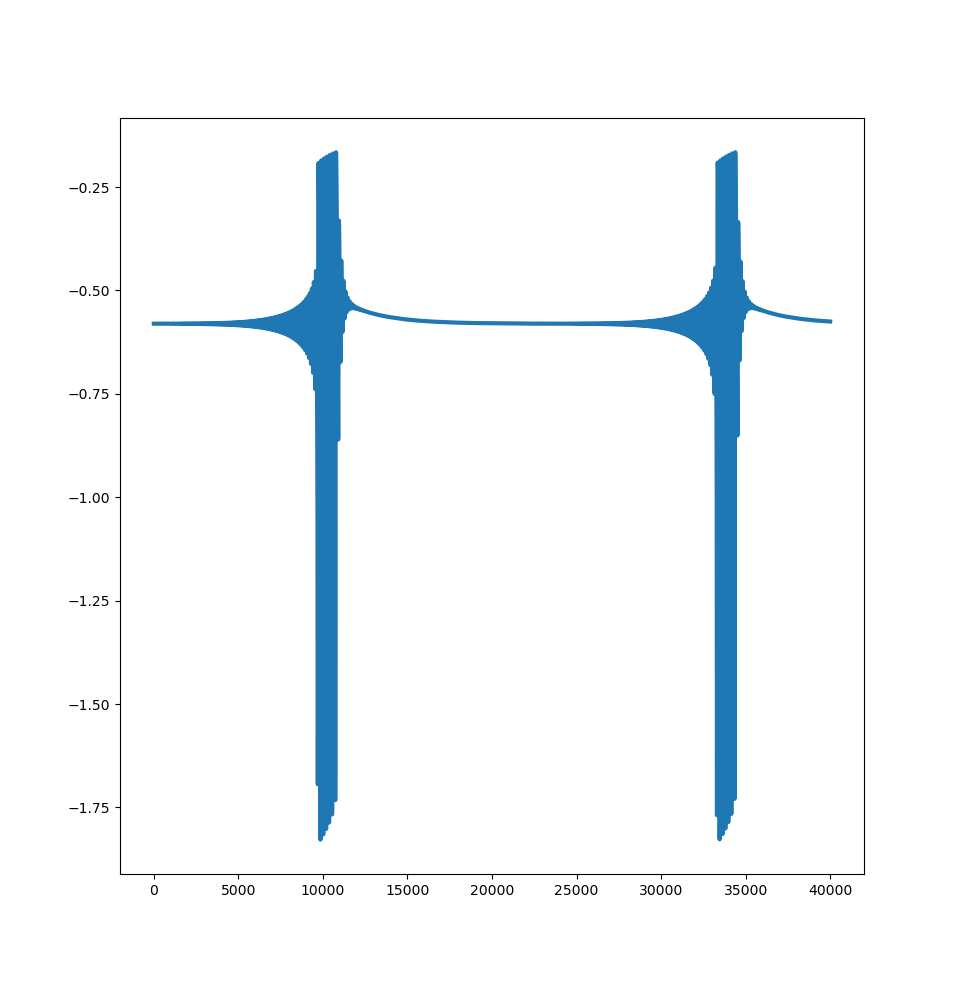}};
	
	\node at (0.15,-2.9) {$t$};
	\node at (-2.9,0) {$x$};
	
	\end{tikzpicture}
		\caption{$c = 1.2$, $k=0.7$, $a=2.2$}
	\end{subfigure}
	\begin{subfigure}[b]{0.48\textwidth}
	\centering
	\begin{tikzpicture}
	\node at (0,0){
	\includegraphics[scale=0.25]{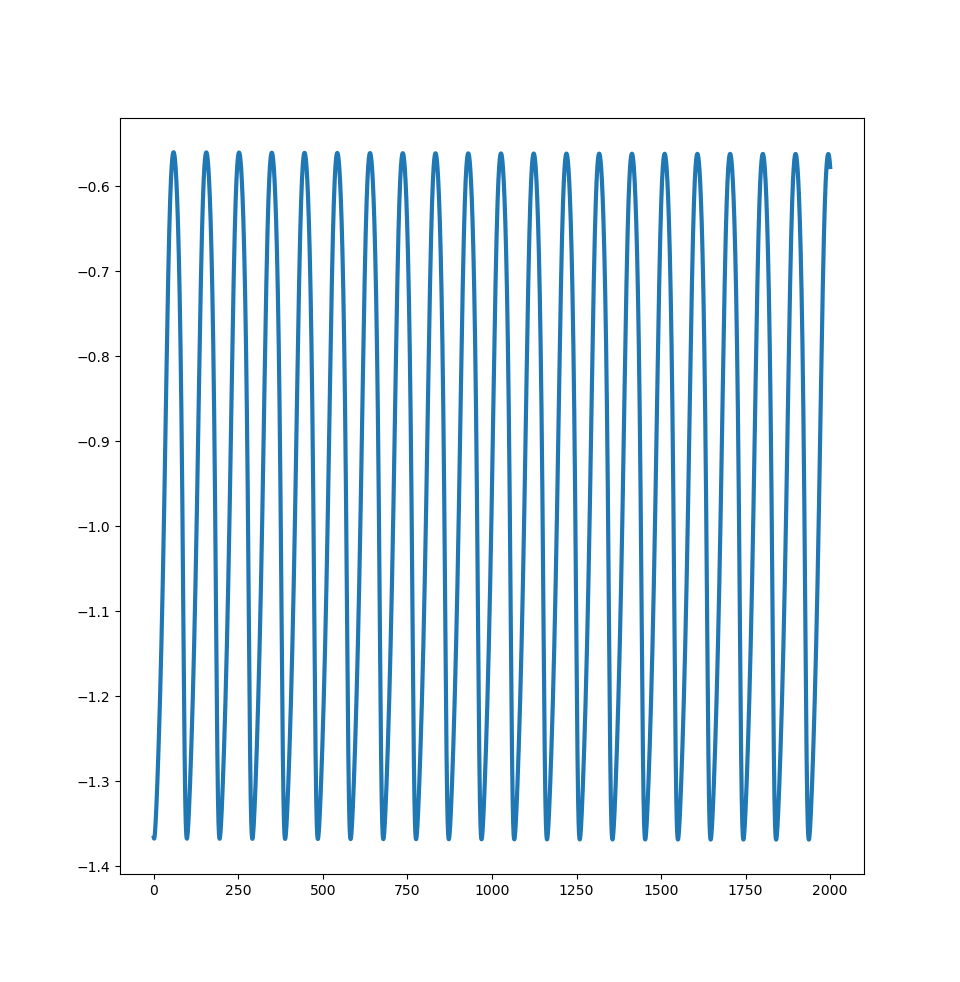}};
	
	\node at (0.15,-2.9) {$t$};
	\node at (-2.9,0) {$x$};
	
	\end{tikzpicture}
	\caption{$c = 1.5$, $k=0.34$, $a=7.78$} 
	\end{subfigure}

	\caption{Examples of possible oscillatory dynamics in Equation~\eqref{elno}: (a) MMO with plateaus, (b) relaxation oscillation with plateaus, (c) MMO with SAOs above, (d) plateau-less relaxation oscillation.}
	\figlab{outline}
\end{figure}
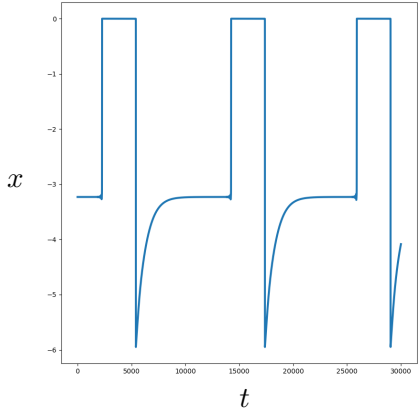
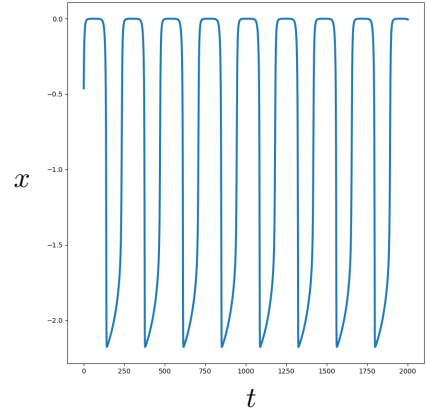
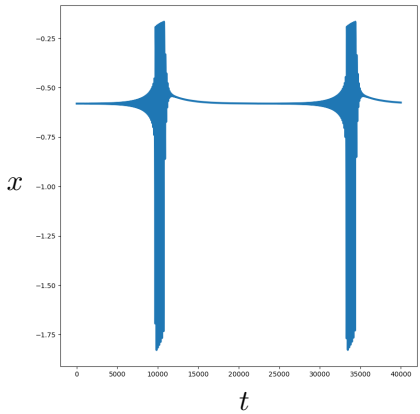
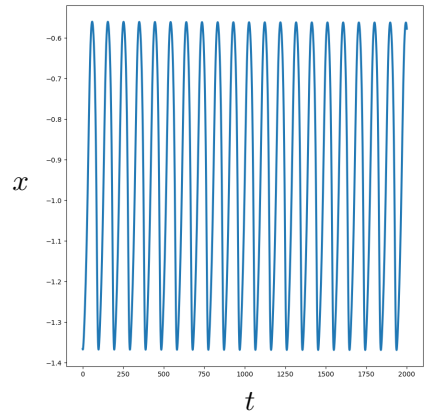

To that end, we identify the dependence of the underlying geometry of \eqref{elno} on the parameters $c$, $k$, and $a$ in the aforementioned three-timescale setting, with $ \delta,\rho >0$ sufficiently small. We first show that the parameter $c$ is associated with the geometric properties of two-dimensional invariant manifolds for \eqref{elno}; following an approach that is similar to the one in \cite{kaklamanos2022bifurcations,kaklamanos2022geometric}, we deduce that the dynamics on these manifolds are described by two dimensional slow-fast systems that can be either in the standard or in the non-standard form of GSPT \cite{wechselberger2020geometric}. For fixed $ c\in (1,c_0) $, the parameter $ k $ is then associated with the geometric properties of one-dimensional immersed invariant submanifolds. The parameter $ a\in(0,a_0) $ is not relevant to the geometry of such invariant (sub)manifolds; rather, it is associated with dynamical phenomena, such as Hopf bifurcations, and with the reduced flow on these (sub)manifolds. In summary, variation in the parameters $c$ and $k$ in \eqref{elno} allows for different types of oscillatory dynamics. The behaviour that is realised from each type then depends on the value of the parameter $a$, cf. \figref{outline}; our aim is, therefore, to associate the various subregions of the $(c,k)$-parameter space to different qualitative behaviours of \eqref{elno} with $\delta, \rho>0$ sufficiently small, as illustrated in \figref{regions}.

\begin{figure}[ht!]
	\centering 
	\begin{tikzpicture}
	\node at (0,0){
	\includegraphics[scale=0.5]{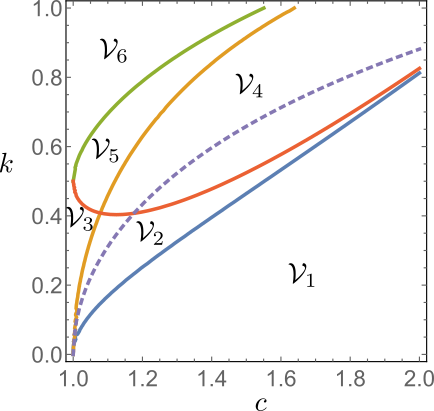}};
	
	\node at (2.6,2.4) {$\textcolor{RoyalPurple}{\mc{C}}$};

	\end{tikzpicture}
	\caption{The $ (c,k) $-plane is divided into six parameter regimes, each one of which corresponds to a type of oscillatory dynamics. Transitions between the possible behaviours in each regime for Equation~\eqref{elno}, with $\delta,\rho>0$ sufficiently small, are controlled by the parameter $a$. The dashed curve $\mc{C}$ distinguishes between two different mechanisms associated with the onset and cessation of oscillatory dynamics in dependence of $a$; see the discussion below \lemmaref{qplus} for details.}
	\figlab{regions}
\end{figure}

The distinction between oscillatory trajectories with different qualitative properties is based on extending the notion of the relative position of folded singularities from \cite{kaklamanos2022bifurcations} to the notion of a relative position of sets where normal hyperbolicity is lost in general. Importantly, following our classification of the underlying singular geometry of \eqref{elno} in the three-timescale regime, we present a wider variety of qualitatively different behaviours in \secref{mains} than are illustrated in \figref{outline}. From our analysis, it will follow that plateau-less relaxation oscillation is equally relevant in the two-timescale context, whereas only relaxation oscillation with plateaus above was documented in \cite{roberts2016mixed}. Finally, we address bifurcation delay in the context of oscillatory dynamics with plateaus above, which was left as an open question in \cite{roberts2016mixed}; specifically, we show that these plateaus arise through delayed passage past an invariant plane in Equation~\eqref{elno}.

Notably, our analysis in this work is complementary to how mixed-mode dynamics in  three-timescale systems are typically analysed, such as in \cite{kaklamanos2022geometric}; that is, instead of explaining previously documented behaviours though the lens of GSPT, we investigate possible geometric configurations of our system, demonstrating the various associated oscillatory behaviours that these configurations can generate. Our approach hence adds value from a wider dynamical systems point of view, as the geometry of Equation~\eqref{elno} extends beyond simple qualitative models for ENSO. As will become apparent, the underlying singular geometry of \eqref{elno} is very similar to the unfolding of Rosenzweig--MacArthur-type systems \cite{duncan2019fast,poggiale2020analysis} and ecosystem models \cite{sadhu2021complex} in three dimensions, as well as to three-dimensional predator-prey models of Lotka-Volterra type \cite{holling1959components,tanner1975stability}. By elucidating the geometry, and hence classifying the dynamics, of \eqref{elno} in the three-timescale regime, we therefore contribute to the understanding of a wider class of multiple-scale models that are ubiquitous in the fields of population dynamics and ecological modelling. 

\medskip

The paper is organised as follows. In \secref{singgeomElno}, we study the singular geometry of Equation~\eqref{elno}, and we investigate the various oscillatory behaviours that can emerge in dependence of the parameters in the system. In \secref{mains}, we relate the resulting singular geometries to the qualitative dynamics of the perturbed system, Equation~\eqref{elno}, with $ \delta,\rho>0 $ sufficiently small. We conclude the paper in \secref{elno-sum} with a brief summary and discussion of our results. 

\section{Singular geometry}
\seclab{singgeomElno}

In this section, we classify the singular geometry of Equation~\eqref{elno} in dependence of the parameters $c$ and $k$, and we demonstrate how the diagram in \figref{regions} is constructed: each region $\mc{V}_i$ in that diagram corresponds to a different type of oscillatory dynamics for $\delta, \rho>0$ sufficiently small in \eqref{elno}. For fixed $c$ and $k$, a particular dynamical scenario out of the possible ones in the corresponding region $\mc{V}_i$ is then realised upon variation of the parameter $a$.

\subsection{The critical manifold $ \mc{M}_1 = \mc{M}_\mc{P}\cup \mc{M}_\mc{S} $}
Given $ \delta>0 $, Equation~\eqref{elno} is written in the fast formulation of GSPT, with the prime denoting differentiation with respect to the fast time $ t $. In the intermediate formulation, i.e. after a rescaling of time as $ \tau = \delta t  $, \eqref{elno} reads 
\begin{subequations}\eqlab{inter}
	\begin{align}
	\delta \dot{x}&=x\big[ x+y+c(1-\tanh{(x+z)})\big]+\rho\delta(x^2-ax), \\
	\dot{y}&=-\rho(ay+x^2), \\
	\dot{z}&=k-z-\tfrac{x}{2};
	\end{align}
\end{subequations}
the slow formulation, after a rescaling of time as $ s = \rho \tau  $ in \eqref{inter}, is given by
\begin{subequations}\eqlab{slow}
	\begin{align}
	\delta\rho \dot{x}&=x\big[ x+y+c(1-\tanh{(x+z)})\big]+\rho\delta(x^2-ax), \\
	\dot{y}&=-ay-x^2, \\
	\rho\dot{z}&=k-z-\tfrac{x}{2}.
	\end{align}
\end{subequations}
In the singular limit of $ \delta =0 $, the so-called layer problem is obtained from Equation~\eqref{elno} as
\begin{subequations}\eqlab{layer}
	\begin{align}
	x'&=x\big[ x+y+c(1-\tanh{(x+z)})\big] =:F(x,y,z), \eqlab{layer-a} \\
	y'&=0, \\
	z'&=0,
	\end{align}
\end{subequations}
while the reduced problem is found by setting $ \delta =0 $ in Equation~\eqref{inter}:
\begin{subequations}\eqlab{reduced}
	\begin{align}
	0&=x\big[ x+y+c(1-\tanh{(x+z)})\big], \eqlab{reduced-a} \\
	\dot{y}&=-\rho(ay+x^2), \\
	\dot{z}&=k-z-\tfrac{x}{2}.
	\end{align}
\end{subequations}
Equilibrium solutions of the one-dimensional layer problem, Equation~\eqref{layer}, define the critical manifold $ \mc{M}_1 =  \mc{M}_\mc{P}\cup \mc{M}_\mc{S}$, where
\begin{align}
\mc{M}_\mc{P} &= \lb (x,y,z)\in \mb{R}^3 ~\lvert~ x = 0\rb \quad\text{and} \eqlab{Mp}\\
\mc{M}_\mc{S} &= \lb (x,y,z)\in \mb{R}^3 ~\lvert~  x+y+c(1-\tanh{(x+z)}) = 0\rb. \eqlab{Ms}
\end{align}
The stability of $ \mc{M}_1 $ is determined by linearisation with respect to the fast variable $ x $ in \eqref{layer-a}:
\begin{align}
F_x= x\lp 1-c\,\tn{sech}^{2}{(x+z)}\rp+\big[ x+y+c(1-\tanh{(x+z)})\big]. \eqlab{xlin}
\end{align} 

By \eqref{Mp}, $ \mc{M}_\mc{P} $ corresponds to the hyperplane $\{ x=0 \}$; the normally hyperbolic subset of $ \mc{M}_\mc{P} $ is therefore defined as 
\begin{align}
\mc{P} = \lb (x,y,z)\in \mc{M}_\mc{P}~\lvert~ F_x\lvert_{x=0}\neq0\rb, 
\end{align}
where, by \eqref{xlin},
\begin{align}
F_x\lvert_{x=0} = y+c(1-\tanh{(z)}) . \eqlab{FMp}
\end{align}
The attracting and repelling subsets $ \mc{P}^a $ and  $ \mc{P}^r $ of $ \mc{P}$ are therefore given by
\begin{subequations}\eqlab{normplanes}
\begin{align}
\mc{P}^{a} &= \lb (x,y,z)\in \mc{M}_\mc{P}~\lvert~ y+c(1-\tanh{(z)})<0\rb\quad\text{and}\eqlab{normplanes-a} \\
\mc{P}^{r} &= \lb (x,y,z)\in \mc{M}_\mc{P}~\lvert~ y+c(1-\tanh{(z)})>0\rb, \eqlab{normplanes-b}
\end{align}
\end{subequations}
respectively. The manifold $ \mc{M}_\mc{P} $ is not normally hyperbolic at $ \mc{F}_{\mc{P}} = {\mc{M}_\mc{P}}\backslash\mc{P}$, where
\begin{align}
\mc{F}_{\mc{P}} &=\lb (x,y,z)\in\mc{M}_1  ~\lvert~ x=0=y+c(1-\tanh{(z)}) \rb; \eqlab{nhMp}
\end{align}
here, we note that $\mc{M}_{\mc{P}}$ and $\mc{M}_{\mc{S}}$ intersect in $\mc{F}_{\mc{P}}$. 

\begin{remark}
	We emphasise that the hyperplane $ \mc{M}_\mc{P} $ is invariant for the full system in \eqref{elno}; moreover, we recall that we will restrict to $ x\leq 0 $ in the following. 
\end{remark}
Similarly, by \eqref{Ms}, the submanifold $ \mc{M}_\mc{S} $ can be written as a graph of $ y $ over $ x $ and $ z $:
\begin{align}
y=  -x-c(1-\tanh{(x+z)}) =:h(x,z) \eqlab{hMs}.
\end{align}
The normally hyperbolic subset of $ \mc{M}_\mc{S} $ is therefore defined as 
\begin{align}
\mc{S} = \lb (x,y,z)\in \mc{M}_\mc{S}~\lvert~ F_x\lvert_{y=h(x,z)}\neq0\rb, 
\end{align}
where, from \eqref{xlin} and \eqref{hMs},
\begin{align}
F_x\lvert_{y=h(x,z)} = x\lp 1-c\,\tn{sech}^{2}{(x+z)}\rp. \eqlab{xFx}
\end{align}
For $ x<0 $, the attracting and repelling subsets $ \mc{S}^a $ and $ \mc{S}^r $ of $\mc{S}$ are therefore given by
\begin{subequations}\eqlab{Sar}
\begin{align}
\mc{S}^a &= \lb (x,y,z)\in \mc{M}_\mc{S}~\lvert~  1-c\,\tn{sech}^{2}{(x+z)}>0\rb\quad\text{and} \eqlab{Sar-a}\\
\mc{S}^r &= \lb (x,y,z)\in \mc{M}_\mc{S}~\lvert~ 1-c\,\tn{sech}^{2}{(x+z)}<0\rb, \eqlab{Sar-b}
\end{align}
\end{subequations}
respectively. (Note that we have accounted for the factor of $x$ in \eqref{xFx}, which is not included in \eqref{Sar}.) The manifold $ \mc{M}_\mc{S} $ loses normal hyperbolicity at $ \mc{F}_{\mc{S}} \cup \mc{F}_{\mc{P}}= {\mc{M}_\mc{S}}\backslash\mc{S}$, where
\begin{align}
\mc{F}_{\mc{S}}&=\lb (x,y,z)\in\mc{M}_1  ~\lvert~ 1-{c}\tn{sech}^{2}{(x+z)} = 0\rb = \mc{L}^-\cup\mc{L}^+; \eqlab{nhMs}
\end{align}
the fold lines $ \mc{L}^\mp $, which separate the normally hyperbolic part of $\mathcal{M}_{\mathcal{S}}$ into $\mc{S}=\mc{S}^{a^-} \cup \mc{S}^{r}\cup \mc{S}^{a^+}$, are given as graphs in the $ (x,z) $-plane  by
\begin{align}
\mc{L}^\mp :\ z = -x\mp\tn{arcsech}\lb \frac{1}{\sqrt{c}}\rb.  \eqlab{Lmp}
\end{align}

The above geometric objects are illustrated in \figref{ms-elno}. In the following, we study the reduced flows on $ \mc{M}_\mc{P} $ and $ \mc{M}_\mc{S} $ under the assumption that these flows are slow-fast systems themselves; in other words, we consider $ \rho>0 $ sufficiently small. 
\begin{figure}[ht!]
	\centering
\includegraphics[scale=0.4]{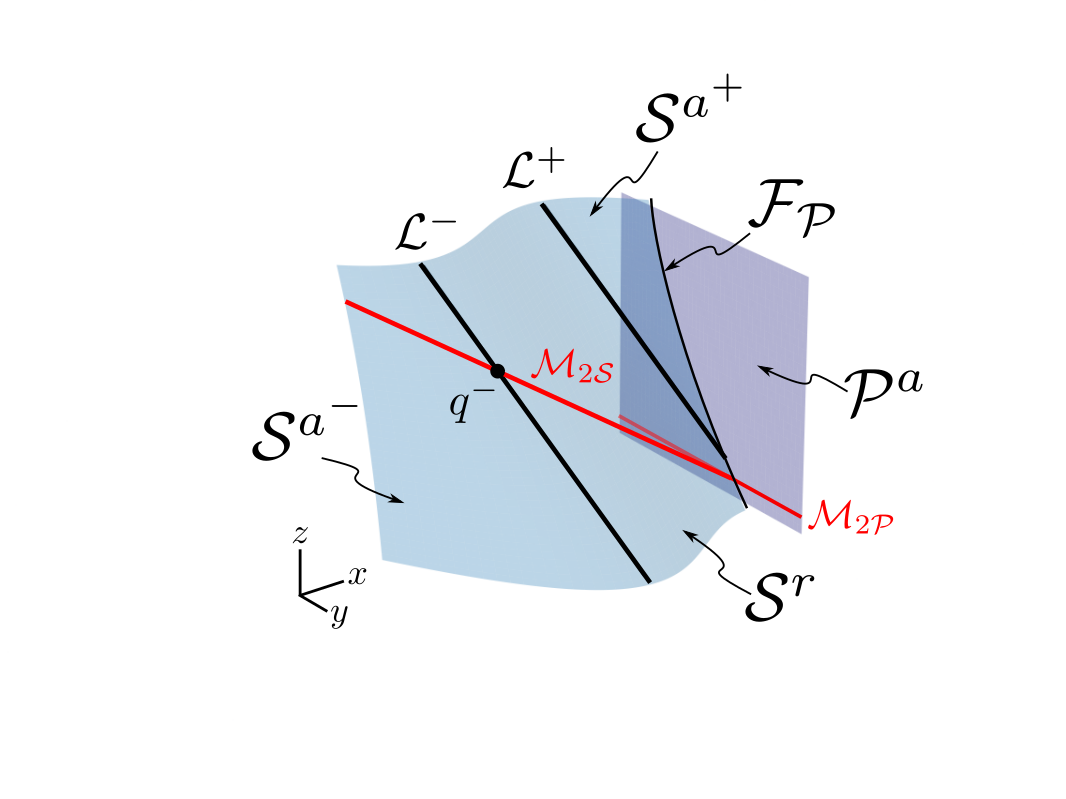}
\caption{Critical manifolds (in grey) and $2$-critical manifolds (in red) for Equation~\eqref{elno}.}
\figlab{ms-elno}
\end{figure}

\subsection{The reduced flow on $ \mc{M}_\mc{P} $}
The reduced flow on $ \mc{M}_\mc{P}	$ is obtained by setting $ x=0 $ in \eqref{reduced}:
\begin{subequations}\eqlab{fastP}
	\begin{align}
	y' &=-\rho ay,\\
	z' &=k-z;
	\end{align}
\end{subequations}
Equation~\eqref{fastP} is linear, with explicit solution
\begin{align}
y(t) = y_0 e^{-\rho a t}\quad \text{and}\quad z(t)=k+(z_0-k)e^{-t} \eqlab{explicit},
\end{align}
from which we can parametrise $y$ with $z$ via
\begin{align}
y(z) = y_0 \lp \frac{z-k}{z_0-k}\rp ^{\rho a}. \eqlab{yofz}
\end{align}
The system in \eqref{fastP} admits a stable node at $ (0,k) $, with eigenvalues $ -1 $ and $ -\rho a $. 

For $ \rho>0 $ sufficiently small, Equation~\eqref{fastP} is a slow-fast system in the standard form of GSPT, where $ y $ is the slow variable and $ z $ is the fast one. In the slow formulation, i.e. after a rescaling of time as $ s = \rho \tau $, \eqref{fastP} reads
\begin{subequations}\eqlab{slowMp}
	\begin{align}
	\dot{y}&=-ay, \\
	\rho\dot{z}&=k-z.
	\end{align}
\end{subequations}
The layer problem on $ \mc{M}_\mc{P} $ is found by setting $ \rho=0 $ in \eqref{fastP}:
\begin{subequations}\eqlab{layP}
	\begin{align}
	y' &=0,\\
	z' &=k-z; \eqlab{layP-b}
	\end{align}
\end{subequations}
its solutions are given by straight fibres with $ y $ constant, as can also be seen by taking $ \rho=0 $ in \eqref{explicit}. Equilibria of \eqref{layP} define the $2$-critical manifold
\begin{align}
\mc{M}_{2\mc{P}} = \lb (x,y,z)\in\mb{R}^3 ~\lvert~ x=0=k-z \rb. \eqlab{M2P}
\end{align}
The reduced flow on $ \mc{M}_{2\mc{P}} $ is obtained by setting $ \rho=0 $ in \eqref{slowMp}:
\begin{subequations}\eqlab{redMp}
	\begin{align}
	\dot{y}&=-ay, \\
	0&=k-z.
	\end{align}
\end{subequations}
Since the fast flow in \eqref{layP-b} is linear with respect to $ z $, the manifold $ \mc{M}_{2\mc{P}} $ is normally hyperbolic -- and, in fact, normally attracting -- everywhere. 

We make the following observation:
\begin{lemma}
	The stable node $(0,0,k)$ of \eqref{elno} -- which corresponds to the point $ (0,k) $ in the context of \eqref{fastP} -- lies on $ \mc{P}^r $ for all $ c,k,a>0 $. 
\end{lemma} 
\begin{proof}
	From \eqref{FMp}, we have that $ F_x\lvert_{(x,y,z)=(0, 0, k)} >0$ for all $ c,k,a>0 $; by \eqref{normplanes}, it then follows that the point $ (0,0,k) $ lies on $ \mc{P}^r $. 
\end{proof}


\subsection{The reduced flow on $ \mc{M}_\mc{S} $}
Next, we describe the reduced flow on the portion $\mc{M}_{\mc{S}}$ of $\mc{M}_1$.
Differentiating the graph representation of $ \mc{M}_1 $ in \eqref{hMs} gives $ \dot{y} = h_x\dot{x}+h_z\dot{z} $; substituting into the reduced problem in \eqref{reduced} and rearranging yields 
\begin{subequations}
	\begin{align*}
	\lp 1-c\tn{sech}^2\lp x+z\rp\rp	\dot{x} &= c\lp k-z-\frac{x}{2}\rp\tn{sech}^2\lp x+z\rp+\rho \lp ah(x,z)+x^2\rp, \\
	\dot{z} &=  k-z-\frac{x}{2},
	\end{align*}
\end{subequations}
which is singular along $ \mc{F}_{\mc{S}} $. Rescaling time by a factor of $ 1-c\,\tn{sech}^2\lp x+z\rp$ in the above gives
\begin{subequations}\eqlab{reducedS}
	\begin{align}
	\dot{x} &= c\lp k-z-\frac{x}{2}\rp\tn{sech}^2\lp x+z\rp+\rho \lp ah(x,z)+x^2\rp, \\
	\dot{z} &= \lp 1-c\,\tn{sech}^2\lp x+z\rp\rp\lp k-z-\frac{x}{2}\rp.
	\end{align}
\end{subequations}
Equation~\eqref{reducedS} is a slow-fast system in the non-standard form of GSPT \cite{wechselberger2020geometric},
\begin{align*}
\begin{pmatrix}
\dot{x} \\ \dot{z}
\end{pmatrix}
=
N(x,z)f(x,z)+\rho\, G(x,z),
\end{align*}
where
\begin{gather*}
N(x,z) = \begin{pmatrix}
c\,\tn{sech}^2\lp x+z\rp\\
1-c\tn{sech}^2\lp x+z\rp
\end{pmatrix}, 
\quad 
f(x,z) = k-z-\frac{x}{2}, 
\quad\text{and}\quad
G(x,z) = \begin{pmatrix}
ah(x,z)+x^2 \\ 0
\end{pmatrix}.
\end{gather*}
The intermediate problem in the double singular limit of $ \rho =0=\delta$ is obtained by setting $ \rho =0$ in \eqref{reducedS}:
\begin{subequations}\eqlab{reducedSsing}
	\begin{align}
	\dot{x} &= c\lp k-z-\frac{x}{2}\rp\tn{sech}^2\lp x+z\rp, \\
	\dot{z} &= \lp 1-c\tn{sech}^2\lp x+z\rp\rp\lp k-z-\frac{x}{2}\rp.
	\end{align}
\end{subequations}

\begin{remark}
In \eqref{reducedS}, we have eliminated the slow variable $ y $, in contrast to \cite{kaklamanos2022bifurcations,kaklamanos2022geometric}, where the intermediate variable was eliminated in the corresponding reduced problems. Alternatively, from the algebraic constraint in \eqref{Ms}, one may choose to express $\mc{M}_{\mc{S}}$ as a graph over $x$ and $y$, with
\begin{align}
z = \tn{arctanh}\lb \frac{x+y}{c}+1\rb - x =:g(x,y). \eqlab{gMs}
\end{align}
Differentiating the representation of $ \mc{M}_\mc{S} $ in \eqref{gMs} gives $ \dot{z} = g_x\dot{x}+g_y\dot{y} $; substituting into the reduced problem in \eqref{reduced} and rearranging then yields
\begin{subequations}
	\begin{align*}
	\lb \frac{1}{c\ls 1-\lp \frac{x+y}{c}+1\rp^2\rs}-1\rb	\dot{x} &= k-\tn{arctanh}\lb \frac{x+y}{c}+1\rb+\frac{x}{2}+\rho\frac{ay+x^2}{c\ls 1-\lp \frac{x+y}{c}+1\rp^2\rs}, \\
	\dot{y} &=  -\rho \lp ay+x^2\rp,
	\end{align*}
\end{subequations}
which is singular along $ \mc{F}_{\mc{S}} $. Rescaling time by a factor of $ \frac{1}{c\ls 1-\lp \frac{x+y}{c}+1\rp^2\rs}-1$ in the above, we find
\begin{subequations}\eqlab{reducedSz}
	\begin{align}
	\dot{x} &= k-\tn{arctanh}\lb \frac{x+y}{c}+1\rb+\frac{x}{2}+\rho\frac{ay+x^2}{c\ls 1-\lp \frac{x+y}{c}+1\rp^2\rs}, \\
	\dot{y} &=  -\rho \lb \frac{1}{c\ls 1-\lp \frac{x+y}{c}+1\rp^2\rs}-1\rb\lp ay+x^2\rp,
	\end{align}
\end{subequations}
which is a slow-fast system in the {standard} form of GSPT \cite{fenichel1979geometric}. {However, as both the representation in \eqref{gMs} and the flow of \eqref{reducedSz} are defined for $(x+y)/c\in(-2,0)$ only, which we cannot guarantee {\it a~priori}, we will not pursue that standard form further here.} 
\remlab{standard}
\end{remark}

Non-stationary solutions of \eqref{reducedSsing} will be called the \textit{intermediate fibres}; for completeness, we note that these fibres can be represented explicitly, as follows.
\begin{lemma}
	For initial conditions $ (x_0, z_0) \in \mc{S}^a$, the intermediate fibres of \eqref{reducedSsing} are given as graphs $ z= \zeta(x;x_0,z_0) $, where
	\begin{align}
	\zeta(x;x_0,z_0) = -x +\tn{arctanh}\lb \frac{x-x_0+c\tanh\lp{x_0+z_0}\rp}{c}\rb . \eqlab{intfib}
	\end{align}
	\lemmalab{zetagraphs}
\end{lemma}

\begin{proof}
	By \eqref{reducedSsing}, for $ k-z-\frac{x}{2}\neq 0 $ we may write
	\begin{align}
	\frac{\tn{d}z}{\tn{d}x} = \frac{1-c\tn{sech}^2\lp x+z\rp}{c\tn{sech}^2\lp x+z\rp} = \frac{1}{c\tn{sech}^2\lp x+z\rp} - 1 = \frac{\cosh^2\lp x+z \rp}{c} - 1. \eqlab{dzdx}
	\end{align}
	As $ u = x+z $, $ \frac{\tn{d}u}{\tn{d}x} = 1+\frac{\tn{d}z}{\tn{d}x}$ gives
	\begin{align*}
	\frac{\tn{d}u}{\tn{d}x} = \frac{\cosh^2\lp u \rp}{c}.
	\end{align*}
	Separating variables, integrating, and reverting to the original coordinates, we obtain the result. 
\end{proof} 

The $2$-critical manifold $ \mc{M}_{2\mc{S}}\subset\mc{S} $ \cite{kaklamanos2022bifurcations} is defined as the set of equilibria of the intermediate problem, Equation~\eqref{reducedSsing}:
\begin{align}
\mc{M}_{2\mc{S}} &= \lb (x,y,z)\in\mb{R}^3  ~\bigg\lvert~ x+y+c\lp 1-\tanh{\lp\frac{x}{2}+k\rp}\rp = 0 = k-z-\frac{x}{2}\rb. \eqlab{M2s}
\end{align}
The Jacobian of the linearisation of \eqref{reducedSsing} about $ \mc{M}_{2\mc{S}}$ has one trivial eigenvalue $ \lambda_0=0 $ and a nontrivial one that is given by
\begin{align*}
\lambda(x,z) &= \la \nabla f, N\ra=-1 +\frac{c \tn{sech}^2\lp x+z\rp}{2},
\end{align*}
where $ \la \cdot, \cdot\ra:\mb{R}^2\to \mb{R} $ denotes the Euclidean inner product in $ \mb{R}^2 $; see \cite{wechselberger2020geometric} for details. Therefore, $ \mc{M}_{2\mc{S}} $ consists of the normally hyperbolic part
\begin{align}
\mc{Z} = \lb (x,y,z)\in\mc{M}_\mc{P}  ~\lvert~ \lambda(x,z) \neq 0\rb,
\end{align}
which can be written as the union $ \mc{Z} = \mc{Z}^a\cup \mc{Z}^r $, where
\begin{align}
\mc{Z}^a = \lb (x,y,z)\in\mc{M}_\mc{P}  ~\lvert~ \lambda(x,z) < 0\rb \quad\text{and}\quad \mc{Z}^r = \lb (x,y,z)\in\mc{M}_\mc{P}  ~\lvert~ \lambda(x,z) > 0\rb.
\end{align}
Here, the folds $ \mc{F}_{\mc{M}_{2\mc{S}}}$ are defined as
\begin{align}
\mc{F}_{\mc{M}_{2\mc{S}}} = \lb (x,y,z)\in\mc{M}_\mc{P}  ~\lvert~ \lambda(x,z) = 0\rb.
\end{align}

Moreover, the folded singularities of $ \mc{M}_\mc{S} $ are defined as the set $Q  = \mc{M}_{2\mc{S}} \cap \mc{L}^\mp=\lb q^-, q^+\rb$, recall \cite{kaklamanos2022bifurcations}. The coordinates of the singularities $ q^\mp=\lp x_{q^\mp}, y_{q^\mp}, z_{q^\mp}\rp$ are 
\begin{gather}
\begin{gathered}
x_{q^\mp} = -2k\mp2\tn{arcsech} \lb \sqrt{\frac{1}{c}}\rb , \quad y_{q^\mp} = h\lp x_{q^\mp},z_{q^\mp}\rp, \quad\text{and}\quad z_{q^\mp} = 2k\pm\tn{arcsech} \lb \sqrt{\frac{1}{c}}\rb,
\end{gathered} \eqlab{qmp}
\end{gather}
respectively, where the function $h$ is defined as in Equation~\eqref{hMs}; see \figref{ms-elno} for an illustration.
Crucially, the above coordinates depend on the parameters $ c $ and $ k $, although we will suppress that dependence in the notation. 
We remark that, as is again apparent from \figref{ms-elno}, the lines $ \mc{L}^\mp $ intersect with $ \mc{F}_{\mc{P}} $, which follows from \eqref{xlin}; similarly, $ \mc{M}_{2\mc{S}} $ intersects with $ \mc{M}_{2\mc{P}} $, which follows from the fact that $ \mc{M}_2 $ is given by the $ x $- and $ z $-nullclines of \eqref{elno}, see \cite{cardin2017fenichel, kaklamanos2022bifurcations,letson2017analysis}. As will become apparent in the following, bounded orbits of \eqref{elno} with $ \delta, \rho>0 $ sufficiently small typically do not interact with these intersections; hence, such orbits can be studied by employing a standard blow-up methodology \cite{krupa2001extending, krupa2001trans}.
%
%
%
%
%

\subsection{Relative locations of non-hyperbolic sets}
In previous works \cite{kaklamanos2022bifurcations,kaklamanos2022geometric}, it was emphasised that the singular geometry of three-timescale systems -- and, in particular, of sets where normal hyperbolicity in the singular limit is lost -- plays an important role for the dynamics of the perturbed flow for $\delta,\rho>0$ sufficiently small, in that it determines the type of oscillation that is possible. Specifically, it was shown that the positions of the folded singularities $ q^\mp $ relative to each other can distinguish between trajectories with different qualitative properties; these relative positions can be classified as follows:
\begin{definition}[\cite{kaklamanos2022mixed,kaklamanos2022bifurcations, kaklamanos2022geometric}]
The folded singularities $q^-$ and $q^+$ are said to be
\begin{enumerate}
    \item (orbitally) remote if $y_{q^-}>y_{q^+}$;
    \item (orbitally) aligned if $y_{q^-}=y_{q^+}$; and
    \item (orbitally) connected if $y_{q^-}<y_{q^+}$.
\end{enumerate}
\defnlab{relpos}
\end{definition}
In words, the folded singularities $ q^\mp $ are remote if no singular cycles of \eqref{elno} exist that evolve on $\mc{M}_{\mc{S}}$ and that pass through both of these singularities; they are connected if there exists a singular cycle which passes through both singularities and which has slow segments in both $\mc{S}^{a^\mp}$; and they are aligned if there exists a singular cycle that passes through both singularities, but that contains no slow segments. Details can be found in \cite{kaklamanos2022bifurcations, kaklamanos2022geometric}. Regarding the position of $ q^\mp $ relative to each other, we have the following result:
\begin{proposition}
	In the double singular limit of $ \rho=0=\delta $, the folded singularities $ q^\mp $ of \eqref{elno} are orbitally remote for all $ c>1 $. 
\end{proposition}
\begin{proof}
In accordance with \defnref{relpos}, the folded singularities $ q^\mp $ are remote for all $ c>1 $ if and only if $y_{q^-}>y_{q^+}$. From \eqref{hMs} and \eqref{qmp}, it follows that
	\begin{align*}
	y_{q^\mp} = -2k\pm 2\tn{arcsech}\lb \frac{1}{\sqrt{c}}\rb-c\lp 1-\tn{tanh}\lb \mp\tn{arcsech}\lb \frac{1}{\sqrt{c}}\rb\rb\rp;
	\end{align*}
	therefore, 
	\begin{align*}
    y_{q^-}-y_{q^+} &= 4\tn{arcsech}\lb   \frac{1}{\sqrt{c}}\rb + \tn{tanh}\lb -\tn{arcsech}\lb \frac{1}{\sqrt{c}}\rb\rb-\tn{tanh}\lb \tn{arcsech}\lb \frac{1}{\sqrt{c}}\rb\rb.
	\end{align*}
	Elementary calculus shows that the right-hand side in the above expression is a positive and increasing function of $ c $ for $ c>1 $, which proves the claim.
\end{proof}
Extending the notions in \defnref{relpos}, we will focus on the position of $ q^\mp $ relative to $ \mc{F}_{\mc{P}} $; we emphasise that both the former and the latter are sets where normal hyperbolicity is lost, albeit in a different fashion. More accurately, we will focus on the projection of $ q^\mp $ onto $ \mc{M}_\mc{P} $ relative to $ \mc{F}_{\mc{P}} $ under the layer flow of \eqref{layer}, which will allow us to describe transitions between qualitatively distinct types of oscillation.

We begin with the following observation:
\begin{lemma}
It holds that
	\begin{align*}
    x_{q^+} <0,
	\end{align*}
	 where $x_{q^+}$ denotes the $x$-coordinate of $q^+$, as before, if and only if 
	\begin{align}
-k+2\tn{arcsech} \lb \sqrt{\frac{1}{c}}\rb <0. \eqlab{kqp}
	\end{align}
	\lemmalab{qplus}
\end{lemma}
\begin{proof}
Follows immediately from \eqref{qmp}.
\end{proof}

Given \lemmaref{qplus}, we define the curve
\begin{align*}
    \mc{C} = \lb (c,k)\in (1,c_0)\times (0,1)~\bigg\lvert~ -k+2\tn{arcsech} \lb \sqrt{\frac{1}{c}}\rb =0\rb, 
\end{align*}
which is illustrated by the dashed curve in \figref{regions}; see also \figref{DAs}(b) below. Therefore, Equation~\eqref{kqp} is satisfied for $ (c,k) $ above the curve $\mc{C}$, which implies that $ q^+ $ lies on the portion of $ \mc{L}^+ $ that is adjacent to $ \mc{S}^{a^+} $ for $ x<0 $. For $ (c,k) $ below the curve $\mc{C}$ in \figref{regions}, $ q^+ $ lies on the portion of $ \mc{L}^+ $ with $ x>0 $; hence, trajectories cannot interact with it due to the invariant plane $\lb x=0\rb$.

As is convention, we denote by $ P(\cdot) $ the projection of a subset of $ \mc{L}^\mp $ along the layer flow of \eqref{layer} onto a portion of the critical manifold $\mc{M}$ and, specifically, onto the first portion of the critical manifold that the fast fibre which emanates from a given point on $\mc{L}^\mp$ intersects with; cf. panels (a) and (c) of \figref{plat?} for an illustration.

\begin{lemma} 
Define
	\begin{align}
\begin{aligned}
A_{q^-}(c,k) := y_{q^-}+c(1-\tanh{(z_{q^-})}),
\end{aligned}
\eqlab{qminP}
	\end{align}
	where $ y_{q^-} $ and $ z_{q^-} $ are given as in \eqref{qmp}. Then, the following statements hold.
	\begin{enumerate}
	    \item $ P(q^-) \in \mc{P}^{a} $ if and only if $ A_{q^-}(c,k)<0 $;
		\item $ P(q^-) \in \mc{F}_{\mc{P}} $ if and only if $ A_{q^-}(c,k)=0 $;
		\item $ P(q^-) \in \mc{S}^{a^+} $ if and only if $ A_{q^-}(c,k)>0 $.
	\end{enumerate}
	\lemmalab{qminrel}
\end{lemma}
\begin{proof}
    Since the $y$- and $z$-coordinates of $P(q^-)$ are $y_{q^-}$ and $z_{q^-}$, respectively, by \eqref{layer}, the second statement follows from the algebraic constraint in \eqref{nhMp}. The remaining two statements are then obtained by substituting $ y_{q^-}$ and $ z_{q^-}$ into \eqref{normplanes}. 
\end{proof}

\begin{figure}[ht!]
	\centering
	\begin{subfigure}[b]{0.3\textwidth}
		\centering
		\includegraphics[scale=0.3]{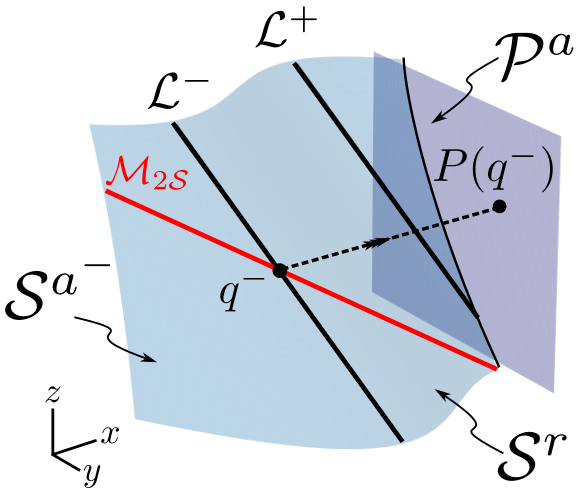}
		\caption{$ A_{q^-}(c,k)<0 $}
	\end{subfigure}
	\begin{subfigure}[b]{0.3\textwidth}
		\centering
		\includegraphics[scale=0.3]{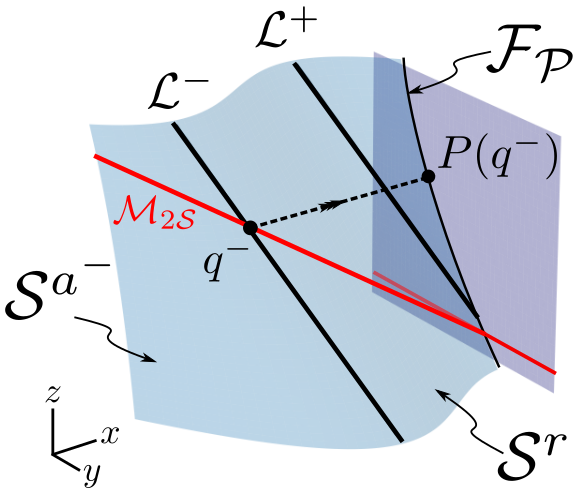}
		\caption{$A_{q^-}(c,k) =0$}
	\end{subfigure}
	\begin{subfigure}[b]{0.3\textwidth}
		\centering
		\includegraphics[scale=0.3]{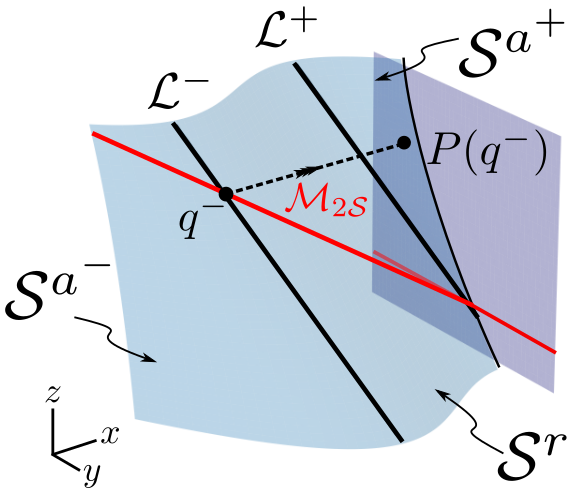}
		\caption{$A_{q^-}(c,k)>0 $}
	\end{subfigure}
	\caption{Illustration of \lemmaref{qminrel}. (a) If $ A_{q^-}(c,k)<0 $, then  $ P(q^-) \in \mc{P}^{a} $, which implies that there exists a singular cycle with a segment evolving in the plane $ \lb x=0\rb $ and not in $ \mc{S}^{a^+} $; in \secref{mains}, the corresponding $ (c,k) $-regime will be associated with oscillatory trajectories that feature plateaus above for $ \delta,\rho>0 $. (c) If $ A_{q^-}(c,k)>0 $, then  $ P(q^-) \in \mc{S}^{a^+} $, which implies that there exists a singular cycle with a segment evolving in $ \mc{S}^{a^+} $ and not in the plane $ \lb x=0\rb $, which will be associated with the existence of oscillatory trajectories without plateaus above, in dependence also of the parameter $ a $. The transition between the two regimes is shown in (b), where $q^-$ is connected to $\mc{F}_{\mc{P}}$ by a fast fibre of \eqref{layer} for $A_{q^-}(c,k)=0$.}
	\figlab{plat?}
\end{figure}

\lemmaref{qminrel} implies that, if the parameters $ c $ and $ k $ satisfy $ A_{q^-}(c,k)=0 $, then the folded singularity $ q^- $ is connected to $ \mc{F}_\mc{P} $ by a fast fibre of \eqref{layer} in the singular limit of $ \delta=0=\rho $, as illustrated in \figref{plat?}(b). More generally, we will denote by $ p_* = (x_*,y_*,z_*)$ the point on $ \mc{L}^- $ that is connected to $ \mc{F}_\mc{P} $ by a fast fibre of \eqref{layer}, i.e. 
\begin{align}
p_* = \lb p \in \mc{L}^- \ \lvert \ P(p) \in \mc{F}_\mc{P}\rb. \eqlab{pstar}
\end{align}
\begin{lemma}
    The point $p_*$, as defined in \eqref{pstar}, exists  and is unique for any fixed $c\in (1,c_0)$.
\end{lemma}
\begin{proof}
    From Equations~\eqref{FMp} and \eqref{Lmp}, on $\mc{M}_\mc{P}$ we obtain
        \begin{align}\eqlab{z}
        \tn{arctanh}\lb \frac{y}{c}+1\rb = z = -\tn{arcsech}\lb \frac{1}{\sqrt{c}}\rb.
    \end{align}
    The left-hand side in \eqref{z} is a monotone function that tends to $\pm\infty$ as $y$ goes to $\pm\infty$, while the right-hand side is constant; therefore, there exists a unique point $(y_*,z_*)$ that solves \eqref{z} for every $c\in (1,c_0)$. Finally, \eqref{Lmp} then implies $x_*=-z_*-\tn{arcsech}\lb \frac{1}{\sqrt{c}}\rb$.
\end{proof}

We note that, since $p_*$ lies on $\mc{L}^-$, it holds that $ x_*<0 $. In addition to the above, we have that
\begin{align*}
P(\mc{L}^-\lvert_{y>y_*}) \subset \mc{S}^{a^+} \quad\text{and}\quad
P(\mc{L}^-\lvert_{y<y_*}) \subset \mc{P}^{a}, 
\end{align*}
which follows from the fact that the left-hand side of \eqref{z} is greater than the right-hand side for $ z<z_* $, and vice versa for $ z>z_* $. Clearly, for $A_{q^-}(c,k)=0$, we have $p_*=q^-$.

Moreover, we will denote by $ q^* = (x^*, y^*, z^*) $ the point on $ \mc{L}^- $ that has the same $y$-coordinate as $ q^+ $, i.e. that lies in the same plane that is parallel to the fast fibres of \eqref{layer} as $ q^+ $:
\begin{align}
q^*=\mc{L}^-\cap\{y=y_{q^+}\}. \eqlab{qstar}
\end{align}
Note that, since the lines $\mc{L}^\mp$ are parallel, and since the point $q^*$ is contained in the plane $\lb y = y_{q^+}\rb$, the latter exists and is unique for all $(c,k)\in (1,c_0)\times(0,1)$.
\begin{figure}[ht!]
	\centering
	\begin{subfigure}[b]{0.3\textwidth}
		\centering
		\includegraphics[scale=0.3]{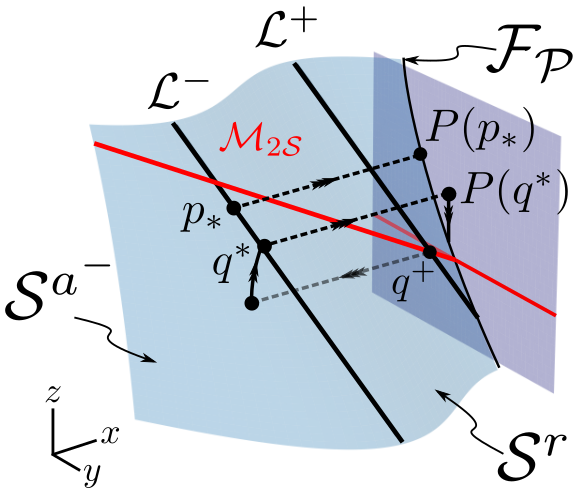}
		\caption{$ A_{q^*}(c,k)>0 $}
	\end{subfigure}
	\begin{subfigure}[b]{0.3\textwidth}
		\centering
		\includegraphics[scale=0.3]{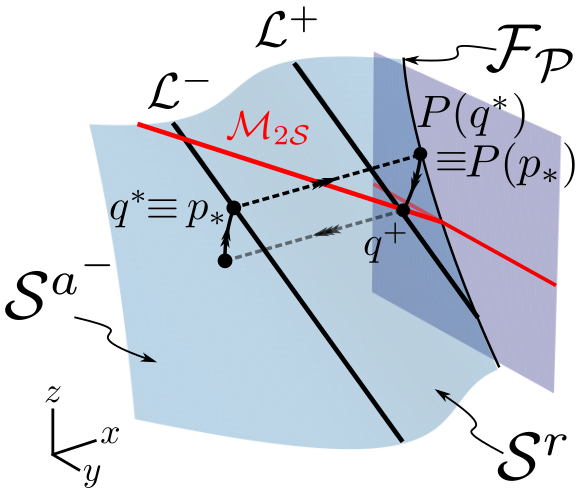}
		\caption{$A_{q^*}(c,k) =0$}
	\end{subfigure}
	\begin{subfigure}[b]{0.3\textwidth}
		\centering
		\includegraphics[scale=0.3]{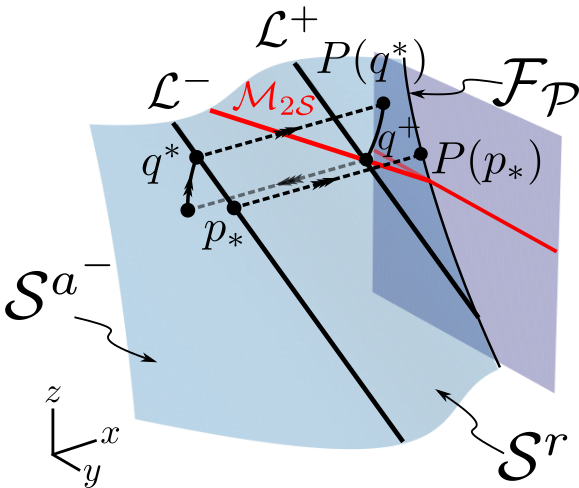}
		\caption{$A_{q^*}(c,k)<0 $}
	\end{subfigure}
	\caption{Illustration of \lemmaref{qstarel}. (a) If $ A_{q^*}(c,k)>0 $, then the location of $ q^+ $ is such that there  exists no singular cycle with endpoint in $ \mc{S}^{a^+} $ that passes through $ q^+ $ -- note that the singular trajectory which emanates from $ q^+ $ does not form a closed orbit. In \secref{mains}, it will be illustrated that MMO trajectories with epochs of perturbed slow dynamics ``above'' in the vicinity of $\mc{M}_{2\mc{S}}$ are not possible in this $ (c,k) $-regime. (c) If $ A_{q^*}(c,k)<0 $, then there exists a singular cycle with endpoint in $ \mc{S}^{a^+} $ that passes through $ q^+ $; in \secref{mains}, it will be illustrated that MMO trajectories with epochs of perturbed slow dynamics ``above'' in the vicinity of $\mc{M}_{2\mc{S}}$ are, in fact, possible in this $ (c,k) $-regime, in dependence of the parameter $ a $. The transition between the two regimes is shown in (b), where $q^*=p^*$ is connected to $ \mc{F}_\mc{P} $ by a fast fibre of \eqref{layer} for $A_{q^*} (c, k) = 0$.}
	\figlab{less?}
\end{figure}

\begin{lemma}
	Define
	\begin{align}
\begin{aligned}
A_{q^*}(c,k) := y_{q^+}-\tn{arcsech}\lb \frac{1}{\sqrt{c}}\rb-\tn{arctanh}\lb \frac{y_{q^+}}{c} +1\rb+c\lp 1 - \tanh\lb -\tn{arcsech}\lb \frac{1}{\sqrt{c}}\rb\rb\rp,
\end{aligned}
	\end{align}
	where $ y_{q^+} $ is as given in \eqref{qmp}. Then, the following statements hold.
	\begin{enumerate}
		\item 	$ P(q^*) \in {\mc{P}^a}$ if and only if $ A_{q^*}(c,k)>0 $;
		\item 	$P(q^*) \in \mc{F}_{\mc{P}} $ if and only if $ A_{q^*}(c,k)=0 $;
		\item 	$  P(q^*) \in \mc{S}^{a^+} $ if and only if $ A_{q^*}(c,k)<0 $.
	\end{enumerate}
	\lemmalab{qstarel}
\end{lemma}
\begin{proof}
	In terms of the second statement, we have that if $P(q^*)\in\mc{F}_{\mc{P}}$, then $y_{q^+} = y_*$ implies
	\begin{align}\eqlab{zast}
	z_* = \tn{arctanh}\lb \frac{y_{q^+}}{c}+1\rb,
	\end{align}
	by the algebraic constraint in \eqref{nhMp}; recall that $ p_* = (x_*, y_*, z_*) $ and $ q^* = (x^*, y^*, z^*) $ are defined in \eqref{pstar} and \eqref{qstar}, respectively. Moreover, 
	\begin{align}\eqlab{xast}
	x_* = -\tn{arcsech}\lb \frac{1}{\sqrt{c}}\rb-z_*, 
	\end{align}
	by the algebraic constraint on $ \mc{L}^- $ in \eqref{Lmp}. Then, we have that
	\begin{align}\eqlab{eq:Azero}
	\begin{split}
	y_{q^+} &= h(x_*, z_*) \\
	&=\tn{arcsech}\lb \frac{1}{\sqrt{c}}\rb+\tn{arctanh}\lb \frac{y_{q^+}}{c} +1\rb-c\lp 1 - \tanh\lb- \tn{arcsech}\lb \frac{1}{\sqrt{c}}\rb\rb\rp,
	\end{split}
	\end{align}
	which, by collecting terms on the left-hand side, gives $ A_{q^*}(c,k)=0 $.

	Conversely, assume that $A_{q^*}(c,k)=0$, i.e. that 
	\eqref{eq:Azero} holds.
	Then, \eqref{Lmp} implies \eqref{xast}, which gives
	\begin{align*}
	y_{q^+} &= -x_*-z_*+\tn{arctanh}\lb \frac{y_{q^+}}{c} +1\rb-c\lp 1 - \tanh\lb x_*+z_*\rb\rp \\
	&=y_{q^+}-z^*+\tn{arctanh}\lb \frac{y_{q^+}}{c} +1\rb,
	\end{align*}
	by \eqref{hMs}. Hence, it necessarily must hold that
	$y_{q^+}+c(1-\tanh(z_*))=0$ which, by \eqref{nhMp}, yields $P(q^*)\in\mc{F}_{\mc{P}}$, as claimed.
	
	The other two statements follow from the properties of $A_{q^*}(c,k)$; see \figref{DAs}(b).
\end{proof}
\lemmaref{qstarel} is illustrated in \figref{less?}. \lemmaref{qminrel} and \lemmaref{qstarel} are summarised in the following corollary; see \figref{DAs}(a) for an illustration.
\begin{cor}
Define
\begin{align*}
\mc{D}_1 &= \lb (c,k)\in (1,c_0)\times (0,1) \ \lvert \ A_{q^-}(c,k)<0 , \ A_{q_*}(c,k)>0 \rb,\\
\mc{D}_2 &= \lb (c,k)\in (1,c_0)\times (0,1) \ \lvert \ A_{q^-}(c,k)>0 , \ A_{q_*}(c,k)>0 \rb,\quad\text{and} \\
\mc{D}_3 &= \lb (c,k)\in (1,c_0)\times (0,1) \ \lvert \ A_{q^-}(c,k)>0 , \ A_{q_*}(c,k)<0 \rb,
 \end{align*}
 as shown in \figref{DAs}. Then, the following statements hold.
 \begin{enumerate}
 	\item If $ (c,k) \in \mc{D}_1$, then $ P(q^-), P(q^*)\in \mc{P}^a $;
 	\item if $ (c,k) \in \mc{D}_2$, then $ P(q^-) \in \mc{S}^{a^+}$, while $  P(q^*)\in \mc{P}^a $;
 	\item if $ (c,k) \in \mc{D}_3$, then $ P(q^-), P(q^*) \in \mc{S}^{a^+}$. 
 \end{enumerate}
\corlab{Ds}
\end{cor}

\corref{Ds} implies that, if $ (c,k) \in \mc{D}_1$, then there exists no singular cycle which passes through $ q^-$ and which has segments that evolve on $ \mc{S}^{a^+} $. Moreover, if $ (c,k) \in \mc{D}_2$, then there exists a singular cycle that passes through $ q^- $ and which has a segment that evolves on $ \mc{S}^{a^+} $, but there exists no singular cycle which passes through $ q^+ $ and has a segment that evolves on $ \mc{S}^{a^+} $. Finally, if $ (c,k) \in \mc{D}_3$, then there exists a singular cycle which passes through $ q^- $ and has a segment that evolves on $ \mc{S}^{a^+} $, as well as a singular cycle which passes through $ q^+ $ and has a segment that evolves on $ \mc{S}^{a^+} $; see \figref{plat?} and \figref{less?}. We will relate the above parameter regimes to properties of oscillatory trajectories in \eqref{elno} for $ \delta,\rho>0 $ sufficiently small in \secref{mains}. We remark that, since the orange curve corresponding to $A_{q_*}(c,k)=0$ in \figref{DAs}(a) always lies above the blue curve corresponding to $A_{q^-}(c,k)=0$ in the $(c,k)$-plane, the combination $\big\{A_{q^-}(c,k)<0 , \ A_{q_*}(c,k)<0\big\}$ is not attainable in \corref{Ds}.

\begin{figure}[ht!]
	\centering
	\begin{subfigure}[b]{0.48\textwidth}
		\centering
		\includegraphics[scale=0.5]{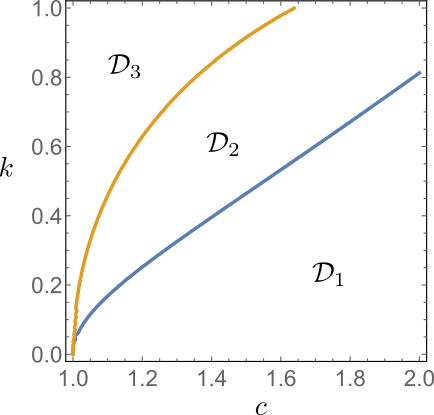}
		\caption{}
	\end{subfigure}
	\begin{subfigure}[b]{0.48\textwidth}
		\centering
		\begin{tikzpicture}
    	\node at (0,0){
    	\includegraphics[scale=0.5]{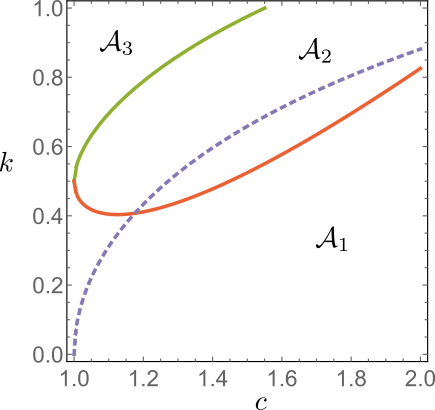}};
    	
    	\node at (2.6,2.4) {$\textcolor{RoyalPurple}{\mc{C}}$};
    
    	\end{tikzpicture}
		
		\caption{}
	\end{subfigure}
	\caption{(a) Parameter regimes described in \corref{Ds}. For $ (c,k)\in \mc{D}_1  $, the projection of the folded singularity $ q^- $ lies in $ \mc{P}^a $; for $ (c,k)\in \mc{D}_2  $, the projection of $ q^- $ lies in $ \mc{S}^{a^+} $, while the projection of the point $q^*$ on $ \mc{L}^- $ with the same $ y $-coordinate as $ q^+ $ lies in $ \mc{P}^a $; for $ (c,k)\in \mc{D}_3  $, the projection of $ q^- $ and of the point $q^*$ lies in $ \mc{S}^{a^+} $, cf. \figref{plat?} and \figref{less?}.\\
		(b) Parameter regimes described in \corref{As}. For $ (c,k)\in \mc{A}_1 $, i.e. below the red curve, there exist $ a^- = a^-(c,k)>0 $ and $ a^+ = a^+(c,k)>0 $, with $ a^->a^+ $, such that an equilibrium point $ \hat{p}_0 $ of the reduced flow, given by \eqref{eqbria}, lies on $ \mc{S}^r $ for $ a\in(a^+,a^-) $. For $ (c,k)\in \mc{A}_2 $, there exists $ a^+ = a^+(c,k)>0 $ such that an equilibrium point of the reduced flow lies on $ \mc{S}^r $ for $ a>a^+ $. For $ (c,k)\in \mc{A}_3  $, there exists no $ a>0 $ such that an equilibrium point of the reduced flow lies on $ \mc{S}^r $. Note that the dashed purple curve $\mc{C}$ is not meant to divide $ \mc{A}_1 $ and $ \mc{A}_2 $ into further subregions; rather, it is related to the location of $ q^+ $, in accordance with \lemmaref{qplus}.}
	\figlab{DAs}
\end{figure}

We now consider the reduced flow on $ \mc{M}_{2\mc{S}} $, which is given by
\begin{align}\eqlab{redS}
\begin{pmatrix}
\dot{x} \\ \dot{z}
\end{pmatrix}
= \left[\frac{\det \lp N\lvert G\rp }{\la \nabla f, N\ra}\begin{pmatrix}
-\partial_z f \\ \partial_x f 
\end{pmatrix}\right];
\end{align}
see \cite{jelbart2020two}. Equilibria on $ \mc{M}_{2\mc{S}} $ are found either by evaluating \eqref{redS}, or by requiring that $ ay+x^2 =0$ in Equation~\eqref{slow}, in addition to the algebraic constraints in \eqref{M2s}, which is equivalent to solving
\begin{align}
a x-x^2+a c \lp 1-\tanh\lp \frac{x}{2}+k\rp\rp=0. \eqlab{eqbria}
\end{align}
We denote the resulting equilibrium point of the reduced flow on $ \mc{M}_{2\mc{S}} $ by 
\begin{align*}
    \hat{p}_0 = (\hat{x}_0 , \hat{y}_0 , \hat{z}_0 ).
\end{align*}
(We note that we have numerically found $ \hat{p}_0 $ to be unique in the parameter regimes considered here.)
It then follows that, in the singular limit of $\rho=0$ in \eqref{slow}, the equilibrium $\hat{p}_0$ lies on $ \mc{S}^{a^-} $ if $ \hat{x}_0<x_{q^-} $, whereas it lies on $ \mc{S}^{r} $ if $ \hat{x}_0>x_{q^-} $, as the lines $\mc{L}^\mp$ are defined by $x$ constant; see e.g.~\figref{plat?}. 

Using the implicit function theorem, one can deduce from \eqref{elno-a} and \eqref{eqbria} that for $ \delta, \rho>0 $ sufficiently small and $ a=\mc{O}(1) $, the point $\hat{p}_0$ lies $ \mc{O}(\delta,\rho) $-close to a true, ``global" equilibrium
\begin{align*}
\hat{p} = (\hat{x}, \hat{y}, \hat{z})
\end{align*}
of the full system, Equation~\eqref{elno}.

\begin{remark}
For $ a\gg1 $, the time-scale separation in the standard form of GSPT in Equation~\eqref{elno} breaks down due to the large $ \mc{O}(a\rho) $-terms on the right-hand side therein. Further investigation of that parameter regime is included in plans for future work.
\remlab{breakdown}
\end{remark}

Solving \eqref{eqbria} for $ a $, we obtain
\begin{align}
a = a(x) = \frac{x^2}{d},\eqlab{ax}
\end{align}
where we have denoted the denominator in the above by
\begin{align}
    d = d(x) = x+c \lp 1-\tanh\lp \frac{x}{2}+k\rp\rp; \eqlab{dx}
\end{align}
we remark that the latter is a decreasing function of $ x $ for $ x\in (x_{q^-}, x_{q^+}) $ and, hence, that $a$ increases with $x$. We further denote
\begin{align}
a^{\mp}(c,k) :=\frac{x_{q^\mp}^2}{d^\mp(c,k)}, \eqlab{amp}
\end{align}
where we have defined
\begin{align}
    d^\mp(c,k) = x_{q^\mp}+c \lp 1-\tanh\lp \frac{x_{q^\mp}}{2}+k\rp\rp. \eqlab{damp}
\end{align}
If the denominator $ d^{-}(c,k) $, respectively $ d^{+}(c,k) $, is positive, then for fixed $ (c,k)\in (1,c_0)\times(0,1) $, the equilibrium point $\hat{p}_0$ is found at $ q^- $, respectively at $ q^+ $, for $ a=a^-(c,k)>0 $, respectively for $ a=a^+(c,k) >0$. (Note that the numerators in \eqref{amp} are always positive.) The graphs of $ d^{\mp}(c,k)=0 $ in the $ (c,k) $-plane are shown in panel (b) of \figref{DAs}. For fixed $(c,k)$, we therefore have $d^{-}(c,k)>d^{+}(c,k)$; we hence distinguish between the following three cases:

\begin{cor}
	Denote 
	\begin{align*}
	\mc{A}_1 &= \lb (c,k)\in (1,c_0)\times (0,1) \ \lvert \ d^{-}(c,k)>0 , \ d^{+}(c,k){>}0 \rb,\\
	\mc{A}_2 &= \lb (c,k)\in (1,c_0)\times (0,1) \ \lvert \ d^{-}(c,k)<0 , \ d^{+}(c,k)>0 \rb,\quad\text{and} \\
	\mc{A}_3 &= \lb (c,k)\in (1,c_0)\times (0,1) \ \lvert \ d^{-}(c,k)<0 , \ d^{+}(c,k)<0 \rb,
	\end{align*}
	as shown in \figref{DAs}. 
	\begin{enumerate}
		\item If $ (c,k) \in \mc{A}_1$, then there exists $ a^\mp= a^\mp(c,k) $, given by $\eqref{amp}$, such that
		\begin{enumerate}
			\item if $ a^+= a^+(c,k) $, then $ \hat{p}_0\equiv q^+ $;
			\item if $ a^-= a^-(c,k) $, then $ \hat{p}_0\equiv q^- $;
			\item if  $ a\in(a^+,a^-) $,  then $\hat{x}_0\in (x_{q^-}, x_{q^+}) $.
		\end{enumerate}
		\item If $ (c,k) \in \mc{A}_2$, then there exists $ a^+= a^+(c,k) $, given by $\eqref{amp}$, such that
		\begin{enumerate}
			\item if $ a^+= a^+(c,k) $, then $ \hat{p}_0\equiv q^+ $;
			\item if  $ a>a^+ $,  then $\hat{x}_0<x_{q^+} $, and there is no $a>0$ such that $ \hat{p}_0\equiv q^-$.
		\end{enumerate}
		\item If $ (c,k) \in \mc{A}_3$, then there is no $a>0$ such that $ \hat{p}_0\in \mc{S}^r $. 
	\end{enumerate}
	\corlab{As}
\end{cor}

(With regard to point 1. in \corref{As}, we note that $a^+(c,k)<a^-(c,k)$ in $\mc{A}_1$ due to $a$ being a strictly increasing function of $x$.)

We now combine the two panels in \figref{DAs} into one figure, dividing the $ (c,k)$-plane into six distinct parameter regimes, as shown in \figref{regions}. 
Note that for $ (c,k) $-values below the dashed purple curve $\mc{C}$ in \figref{DAs}(b), we have $ x_{q^+}>0$ by \lemmaref{qplus}, with $ q^+ $ lying to the right of the plane $ \lb x=0\rb $; therefore, as will become apparent in the following, for $ (c,k) $-values below the curve $\mc{C}$ in \figref{DAs}(b), the value $ a^+ $ is irrelevant: although a stable equilibrium exists for $a\in(0,a^+)$, trajectories with $x\leq0$ cannot reach it.

Finally, we remark that, for fixed $(c,k)\in \mc{V}_2$, there exists an $ a $-value, denoted by $ a_p $, for which the equilibrium $ \hat{p}_0 $ of the reduced flow lies in the plane given by $ \lb y = y_*\rb $; cf.~\eqref{pstar}. We now define the following plane, {which approximates the (invariant) unstable manifold of the global equilibrium $\hat{p}$ to leading order}:
\begin{align*}
\mc{W}(\hat{p}_0) := \lb (x,y,z)\in \mb{R}^3~\lvert~y = \hat{y}_0\rb.
\end{align*}
Then, for $ a<a_p $, it holds that $ P(\mc{L}^-\cap \mc{W}(\hat{p}_0)) \in \mc{P}^a$, whereas $ P(\mc{L}^-\cap \mc{W}(\hat{p}_0)) \in \mc{S}^{a^-}$ for $ a>a_p $, which again follows from the fact that $a$ is an increasing function of $x$ for $x\in(x_{q^-}, x_{q^+})$. In \secref{mains} below, it will become apparent that the $a$-value $a_p$ {distinguishes, in a first approximation}, between oscillatory trajectories that either do or do not feature plateaus above for $(c,k)\in\mc{D}_2$ in \eqref{elno}. {(The transition between the two regimes for $\delta, \rho>0$ sufficiently small will, in fact, be continual rather than abrupt.)}

\section{Outline of dynamics}
\seclab{mains}
In this section, we discuss the perturbed dynamics of Equation~\eqref{elno} for $ \delta,\rho>0 $ and sufficiently small in dependence of the parameters $ c $, $ k $ and $ a $. In particular, we give a qualitative classification of the oscillatory dynamics that will arise upon variation of these parameters in \eqref{elno}.

\subsection{Perturbed dynamics and delayed loss of stability} 

In the following, we consider how the various portions of the critical manifolds $ \mc{M}_{\mc{S}} $ and $ \mc{M}_{\mc{P}} $ perturb for $ \delta, \rho>0 $ sufficiently small. Then, we describe dynamical phenomena, such as delayed loss of stability and Hopf bifurcation, that occur along these perturbed manifolds.

\subsubsection{Normally hyperbolic regime}
By GSPT \cite{cardin2017fenichel,fenichel1979geometric,kaklamanos2022bifurcations,letson2017analysis}, we have that for $ \delta,\rho >0$ sufficiently small, there exist invariant ``slow" manifolds $\mc{S}^{a^\mp}_{\delta\rho}$, $ \mc{S}^r_{\delta\rho}$, $\mc{P}^a_{\delta\rho}$, and $ \mc{P}^r_{\delta\rho}$. The perturbed manifolds $ \mc{S}^{a^\mp}_{\delta\rho} $ and $ \mc{S}^r_{\delta\rho} $ are diffeomorphic, and lie $ \mc{O}(\delta+\rho) $-close in the Hausdorff distance, to their unperturbed, normally hyperbolic counterparts $ \mc{S}^{a^\mp} $ and  $ \mc{S}^r $, respectively. Since $ \mc{M}_{\mc{P}} $ is invariant for Equation~\eqref{elno} for any choice of $ \delta $ and $ \rho $, we conclude that $ \mc{P}^a_{\delta\rho} \equiv \mc{P}^a $ and $ \mc{P}^r_{\delta\rho} \equiv \mc{P}^r $. These manifolds are locally invariant under the flow of \eqref{elno}.

Moreover, for $ \delta,\rho >0$ sufficiently small, there exist invariant ``super-slow" manifolds $ \mc{Z}^a_{\delta\rho} $ and $ \mc{Z}^{r}_{\delta\rho} $ that are diffeomorphic, and $ \mc{O}(\rho) $-close in the Hausdorff distance, to their unperturbed counterparts, the 2-critical manifolds $ \mc{Z}^a $ and $ \mc{Z}^{r} $, respectively. These manifolds are again locally invariant under the flow of \eqref{elno}. (We note that $ \mc{M}_{2\mc{P}}$ is again invariant for any $\delta$ and $\rho$, and that it hence equals its perturbed counterpart.)

\subsubsection{Loss of normal hyperbolicity}
We begin by describing the behaviour of trajectories in the vicinity of $ \mc{F}_\mc{P} $.

In the perturbed system, Equation~\eqref{elno} with $ \delta, \rho>0 $ sufficiently small, the slow sheets $ \mc{S}^r_{\delta\rho}$ and $ \mc{S}^{a^+}_{\delta\rho}$ are ``detached'' 
from the plane $ \lb x=0\rb $, i.e. from $\mc{M}_{\mc{P}}$. That is, although $\mc{M}_{\mc{P}}$ persists for $ \delta, \rho>0 $ sufficiently small, as it is invariant under the flow of \eqref{elno} for all $ \delta$ and $\rho$, the sheets $ \mc{S}^r$ and $ \mc{S}^{a^+}$ perturb to $ \mc{S}^r_{\delta\rho}$ and $ \mc{S}^{a^+}_{\delta\rho}$, respectively, away from $ \mc{F}_{\mc{P}} $. When trajectories on $ \mc{S}^r_{\delta\rho}$ reach the vicinity of $ \mc{F}_{\mc{P}} $, they exit following the fast flow of \eqref{elno}; similarly, trajectories on $ \mc{S}^{a^+}_{\delta\rho}$ which, extended backward in time, reach the vicinity of $ \mc{F}_{\mc{P}} $, again follow the fast flow. Due to the three-timescale structure of Equation~\eqref{elno} for $ \delta, \rho>0 $ sufficiently small, with $ y $ varying slowly, the requisite blow-up transformation near $ \mc{F}_{\mc{P}} $ is analogous to that of a two-dimensional Rosenzweig--MacArthur model which is ``unfolded'' along the $ y $-direction; we refer to \cite[Section 5]{poggiale2020analysis} for details. (We remark that, while the analysis in our case is similar to that of the transcritical singularity in \cite{krupa2001trans}, there are qualitative differences due to the invariance of the plane $\lb x=0\rb$ here, resulting in different weights in the corresponding blow-up transformation.)

Hence, for $ \delta, \rho>0 $ small, trajectories of \eqref{elno} that enter an $ \mc{O}(\delta) $-neighbourhood of $ \mc{P}^a_{\delta\rho} $, but $ \mc{O}(\sqrt{\delta}) $-away from $ \mc{F}_\mc{P} $, follow the slow flow of \eqref{fastP} thereon. {(Here, the restriction away from an $\mc{O}(\sqrt{\delta})$-neighbourhood of $ \mc{F}_\mc{P} $ is due to the rescalings introduced in the local analysis; cf.~again \cite[Section 5]{poggiale2020analysis}.)}
After passing through the vicinity of $ \mc{F}_\mc{P} $,  and instead of being ``immediately'' repelled away from an $ \mc{O}(\delta) $-neighbourhood of $ \mc{P}^r_{\delta\rho} $, trajectories then follow the slow flow of \eqref{fastP} until the accumulated attraction to $ \mc{P}^a_{\delta\rho} $ is balanced by repulsion from $ \mc{P}^r_{\delta\rho} $. Specifically, given an entry point with $z$-coordinate $z_{\rm in}$, the $z$-coordinate $z_{\rm out}$ of the corresponding exit point is calculated using the \textit{way-in/way-out} function:
\begin{align}
\int_{z_{\rm in}}^{z_{\rm out}}\frac{F_x\lvert_{x=0}}{k-z}\tn{d}z = 0, \eqlab{winout}
\end{align}
recall \eqref{FMp} and see \cite{DEMAESSCHALCK20081448,neishtadt2009stability,schecter1985persistent}; this phenomenon, which is also known as \textit{Pontryagin's delay of stability loss}, has been identified in a related system with self-intersecting critical manifold in \cite{sadhu2021complex}. Further, we remark that delayed loss of stability also occurs in the two-dimensional Rosenzweig--MacArthur model studied in \cite{duncan2019fast}, but that it was not addressed there. 

\begin{lemma} 
	Given $ \delta,\rho>0 $ sufficiently small, consider a point $ (x_{\rm in},y_{\rm in}, z_{\rm in})$ in an $ \mc{O}(\delta) $-neighbourhood of $ \mc{P}^a_{\delta\rho} $, but outside an $ \mc{O}(\sqrt{\delta}) $-neighbourhood of $ \mc{F}_\mc{P} $. Then, the trajectory of Equation~\eqref{elno} with initial condition $ (x_{\rm in},y_{\rm in}, z_{\rm in})$ leaves an $ \mc{O}(\delta) $-neighbourhood of $ \mc{P}^a_{\delta\rho} $ at a point $ (x_{\rm out}, y_{\rm out}, z_{\rm out})$ for which
	\begin{align}
	W(z_{\rm in}, z_{\rm out}) := \int_{z_{\rm in}}^{z_{\rm out}}\frac{y_{\rm in} \lp \frac{z-k}{z_{\rm in}-k}\rp ^{\rho a}+c\lp 1-\tanh(z)\rp}{k-z}\tn{d}z = 0 \eqlab{zout}
	\end{align}
	holds. 	
\end{lemma}
\begin{proof}
	The result is based on \cite{DEMAESSCHALCK20081448, schecter1985persistent} and follows from \eqref{winout}, in conjunction with \eqref{FMp} and \eqref{yofz}.
\end{proof}
\noindent We remark that the exit point $(x_{\rm out},y_{\rm out}, z_{\rm out})$ defined by \eqref{zout} can be given in implicit form only in our case. 

We now turn our attention to the behaviour of the perturbed system in \eqref{elno} near $ \mc{F}_{\mc{S}} $ for $ \delta,\rho>0 $ sufficiently small. We focus on the dynamics in a neighbourhood of $ \mc{L}^- $ here; the description of the dynamics near $ \mc{L}^+ $ is similar. 

When trajectories on $ \mc{S}^{a-}_{\delta\rho}$ reach the vicinity of the fold line $ \mc{L}^- $ away from the folded singularity $ q^- $, they ``jump'' to the opposite attracting sheet $ \mc{S}^{a+}_{\delta\rho}$ or $ \mc{P}^a $ following the fast flow; see \cite{wechselberger2005existence,szmolyan2001canards}. On the other hand, when trajectories are attracted to the vicinity of $ q^- $ or to appropriate subregions of $ \mc{Z}^a_{\delta\rho} $, they undergo {epochs of perturbed slow dynamics}. Namely, $ \mc{Z}^a_{\delta\rho} $ can be decomposed into nodally and focally attracting regimes. If trajectories are attracted to the latter regime, then they undergo SAOs of bifurcation delay type, whereas if trajectories are attracted to the former, {then typically no oscillation with discernible amplitude occurs}; see \cite{kaklamanos2022bifurcations} for details. 

\subsubsection{Singular Hopf bifurcation}

Here, we discuss the distinction between steady-state behaviour and oscillatory dynamics in Equation~\eqref{elno}, in dependence of the parameter $ a $, for fixed $ (c,k) \in (1,c_0)\times(0,1)$. We first observe that, near the fold lines $ \mc{L}^\mp $ and away from $ \lb x=0\rb$, Equation~\eqref{elno} can be transformed either into the extended prototypical example studied in \cite{kaklamanos2022bifurcations} or into the canonical form formulated in \cite{letson2017analysis}. 

It then follows that, for $ (c,k)\in \mc{A}_1 $, Equation~\eqref{elno} with $ \delta,\rho>0 $ sufficiently small undergoes singular Hopf bifurcations for $ a = a^-(c,k)+\mc{O}(\delta, \rho) $, since for $ a = a^-(c,k) $, an equilibrium of the reduced flow in \eqref{redS} crosses the fold line $ \mc{L}^- $ in the singular limit of $ \delta=0=\rho $. 

Similarly, for $ (c,k)\in \mc{A}_2 $, for which also \eqref{kqp} holds, i.e. above the dashed purple curve $\mc{C}$ in \figref{DAs}(b), \eqref{elno} undergoes singular Hopf bifurcations for $ a = a^+(c,k)+\mc{O}(\delta, \rho) $ and $ \delta,\rho>0 $ sufficiently small, since for $ a = a^+(c,k) $ and $ \delta=0=\rho $, an equilibrium of the reduced flow in \eqref{redS} crosses the fold line $ \mc{L}^+ $ in the negative $ x$-orthant. 

For $ (c,k)\in \mc{A}_3 $, there exists no $ a>0 $ for which an equilibrium of the reduced flow in \eqref{redS} crosses a fold line $ \mc{L}^\mp $ in the singular limit of $ \delta=0=\rho $ in the negative $ x$-orthant; therefore, the flow of \eqref{elno} converges to steady state for all $ a>0 $. 

\subsection{Oscillatory trajectories}
In this subsection, we present the main qualitative results of this work, summarising the oscillatory dynamics of \eqref{elno} for $ \delta,\rho>0 $ sufficiently small. We combine panels (a) and (b) of \figref{DAs} into \figref{regions}, thus further subdividing the $ (c,k) $-plane, and we illustrate the dynamics of \eqref{elno} for representative $ (c,k) $-values in each of these regimes.

In all numerical simulations below, we consider $ \delta = 0.01=\rho $. 

\subsubsection{\underline{$ (c,k) \in \mc{V}_1 = \mc{D}_1\cap\mc{A}_1$}}
Fix  $ (c,k) \in \mc{V}_1 $, as shown in \figref{regions}.

By \corref{Ds}, it holds that $ P(q^-), P(q^*) \in \mc{P}^a$, i.e. that the projections of both the folded singularity $q^-$ and the associated point $q^\ast$ under the layer flow of \eqref{layer} lie in $\mc{P}^a$. Moreover, by \eqref{amp}, there exist $a$-values $ a^- = a^-(c, k) > 0 $ and $ a^+ = a^+(c, k) > 0 $, respectively, for which the equilibrium $\hat{p}_0 $ given by solving \eqref{eqbria} coincides with the folded singularities $ q^- $ and $ q^+ $, respectively. Correspondingly, for $ a\in(a^+,a^-) $, we have that $ \hat{p}_0\in \mc{S}^r $. However, by \lemmaref{qplus}, it follows that $x_{q^+}>0$ for $(c,k)\in\mc{V}_1$, since this regime lies ``below" the dashed curve $\mc{C}$ in \figref{regions} and \figref{DAs}(b); hence, the Hopf bifurcation that occurs for $a = a^+$ at $q^+$ is irrelevant, as trajectories in the negative $x$-orthant cannot interact with it due to the invariance of $\lb x=0\rb$. Since the reduced flow on $ \mc{M}_{2\mc{S}}\cap\mc{S}^{a^-} $ is directed towards $q^-$ in the absence of a stable equilibrium, singular cycles of Equation~\eqref{elno} contain segments which evolve on $\mc{M}_{\mc{P}}$ for $ a\in(0,a^-) $. For $ \delta,\rho>0$ sufficiently small, \eqref{elno} then exhibits oscillatory dynamics with plateaus above when $ 0<a<a^-+\mc{O}(\delta,\rho) $, with trajectories experiencing delayed loss of stability in their passage along $\mc{M}_{\mc{P}}$. 

In particular, taking $ (c,k) = (1.4,0.2) $ for verification, we calculate that $ a^- \simeq 4.57 $. Numerically, we observe that the flow of Equation~\eqref{elno} converges to steady state for $ a\gtrsim4.5 $, which is consistent with our choice of $\delta$ and $\rho$. Simulated sample trajectories for $ (c,k)\in\mc{V}_1 $ fixed and various values of $ a $, with $ \delta=0.01=\rho $, are illustrated in \figref{V1}.

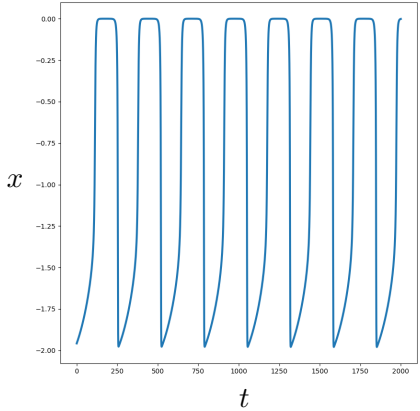
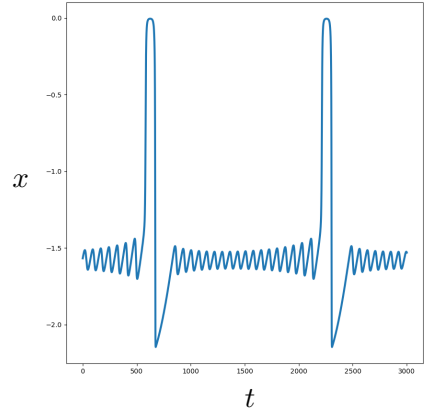
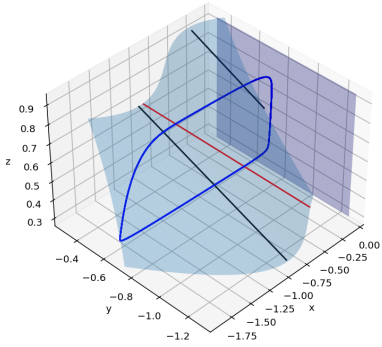
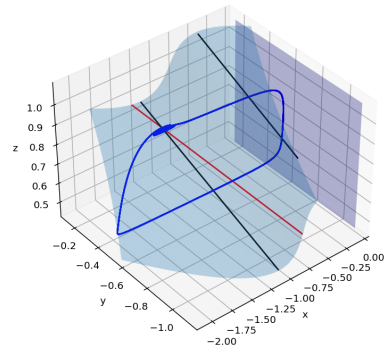
\begin{figure}[ht!]
	\centering
	\begin{subfigure}[b]{0.48\textwidth}
	\centering
	\begin{tikzpicture}
	\node at (0,0){
	\includegraphics[scale=0.25]{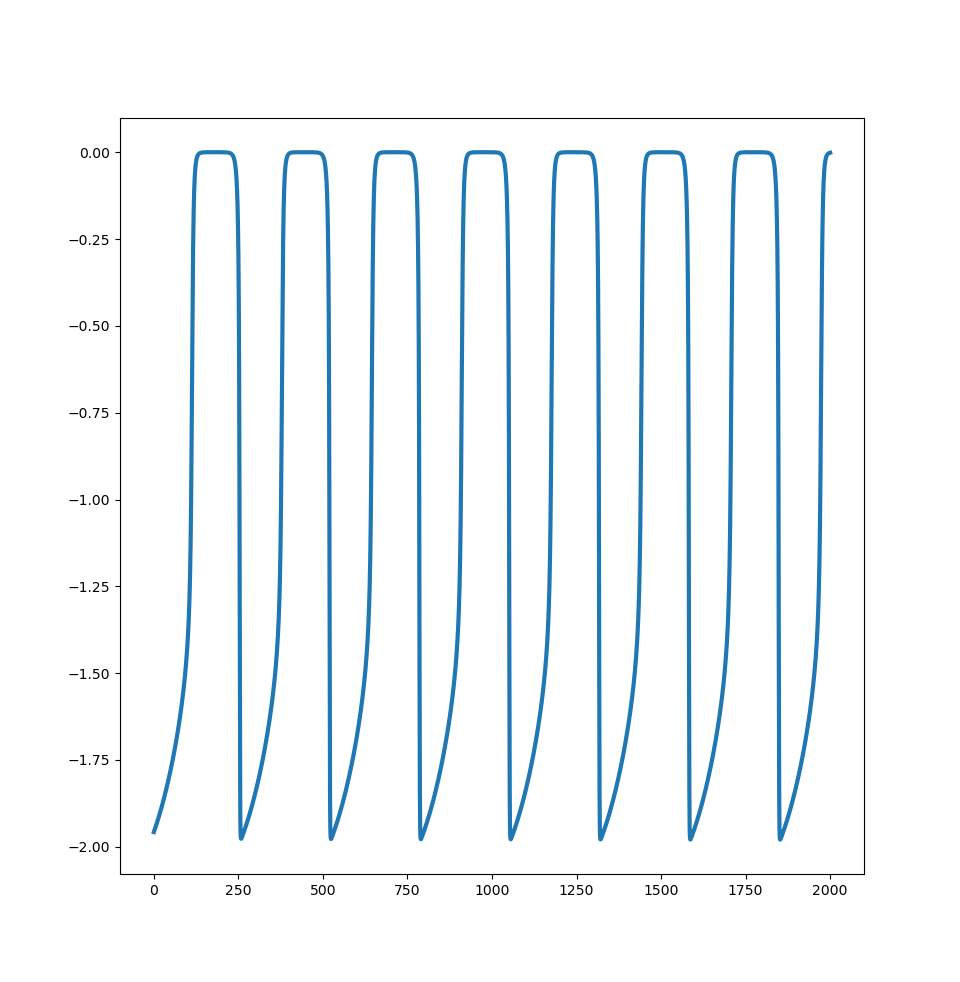}};
	
	\node at (0.15,-2.9) {$t$};
	\node at (-2.9,0) {$x$};
	
	\end{tikzpicture}
		\caption{$c = 1.4$, $k=0.2$, $a=2$}
	\end{subfigure}
	\begin{subfigure}[b]{0.48\textwidth}
	\centering
	\begin{tikzpicture}
	\node at (0,0){
	\includegraphics[scale=0.25]{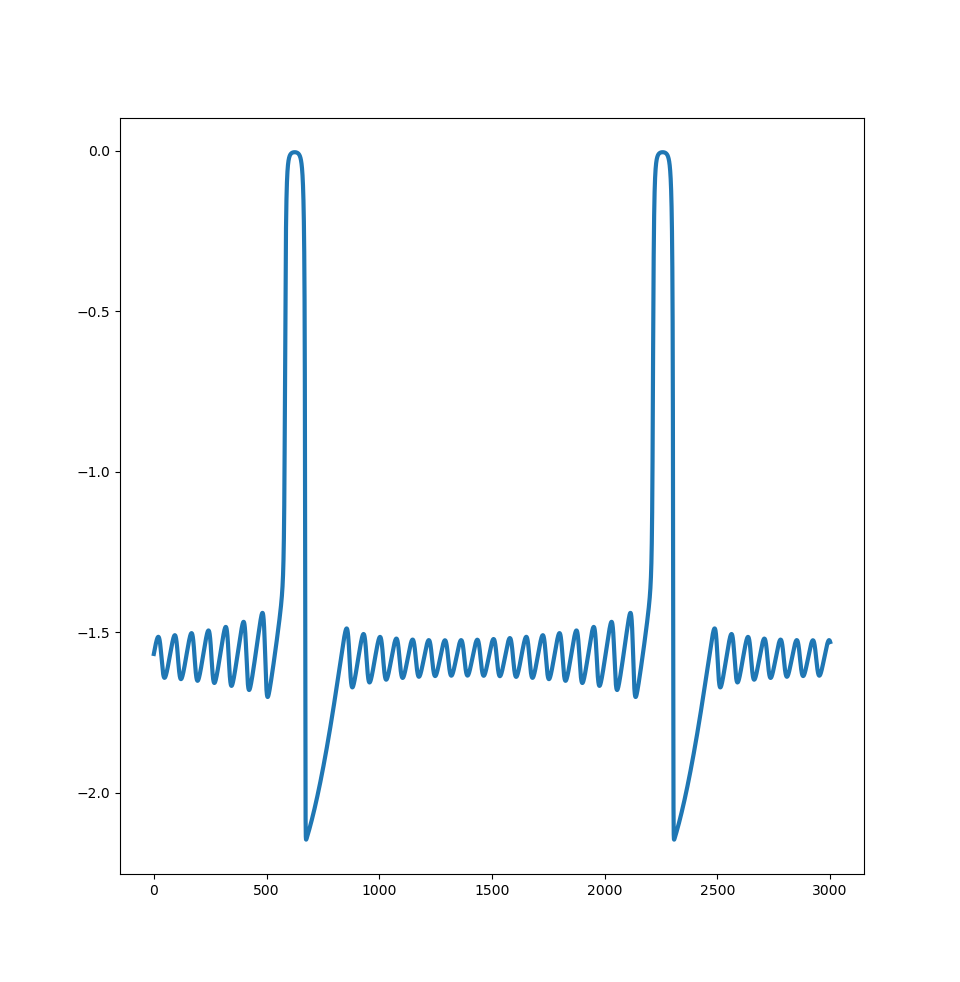}};
	
	\node at (0.15,-2.9) {$t$};
	\node at (-2.9,0) {$x$};
	
	\end{tikzpicture}
		\caption{$c = 1.4$, $k=0.2$, $a=4.4$}
	\end{subfigure}
\\
\begin{subfigure}[b]{0.48\textwidth}
	\centering
	\includegraphics[scale=0.2]{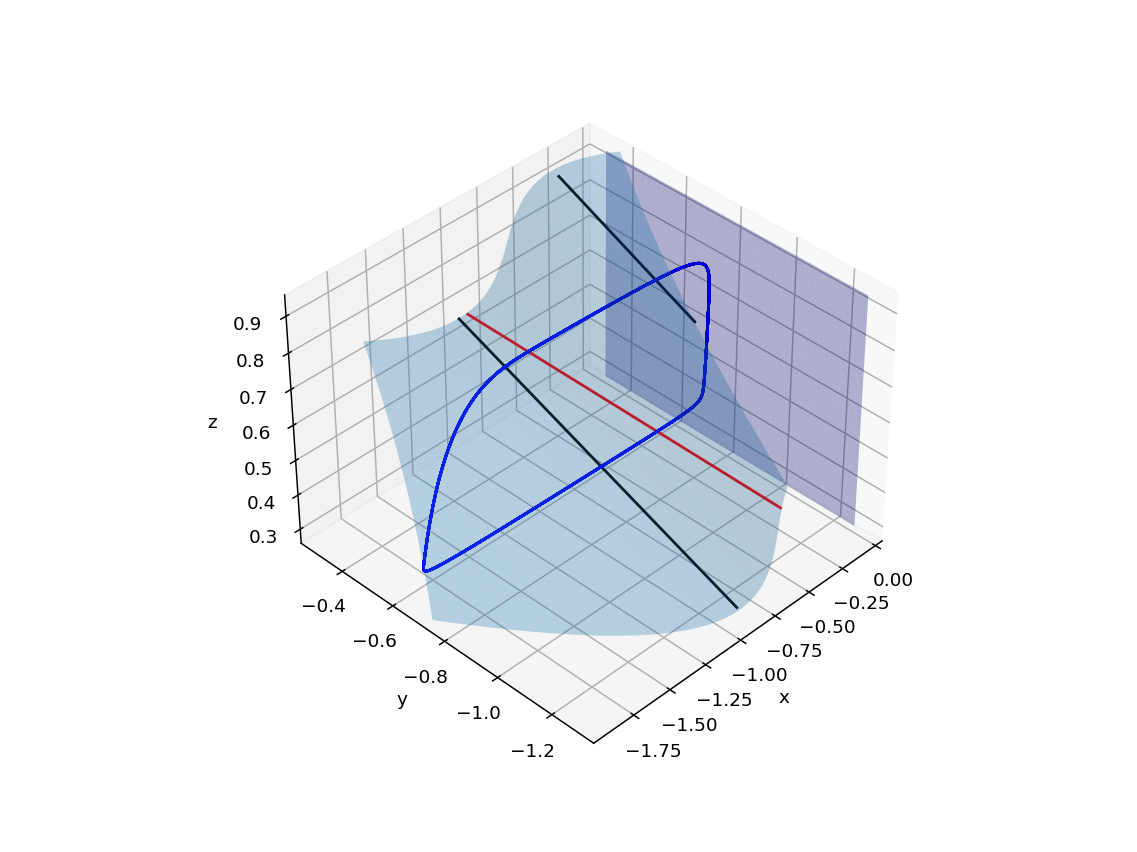}
	\caption{$c = 1.4$, $k=0.2$, $a=2$}
\end{subfigure}
\begin{subfigure}[b]{0.48\textwidth}
	\centering
	\includegraphics[scale=0.2]{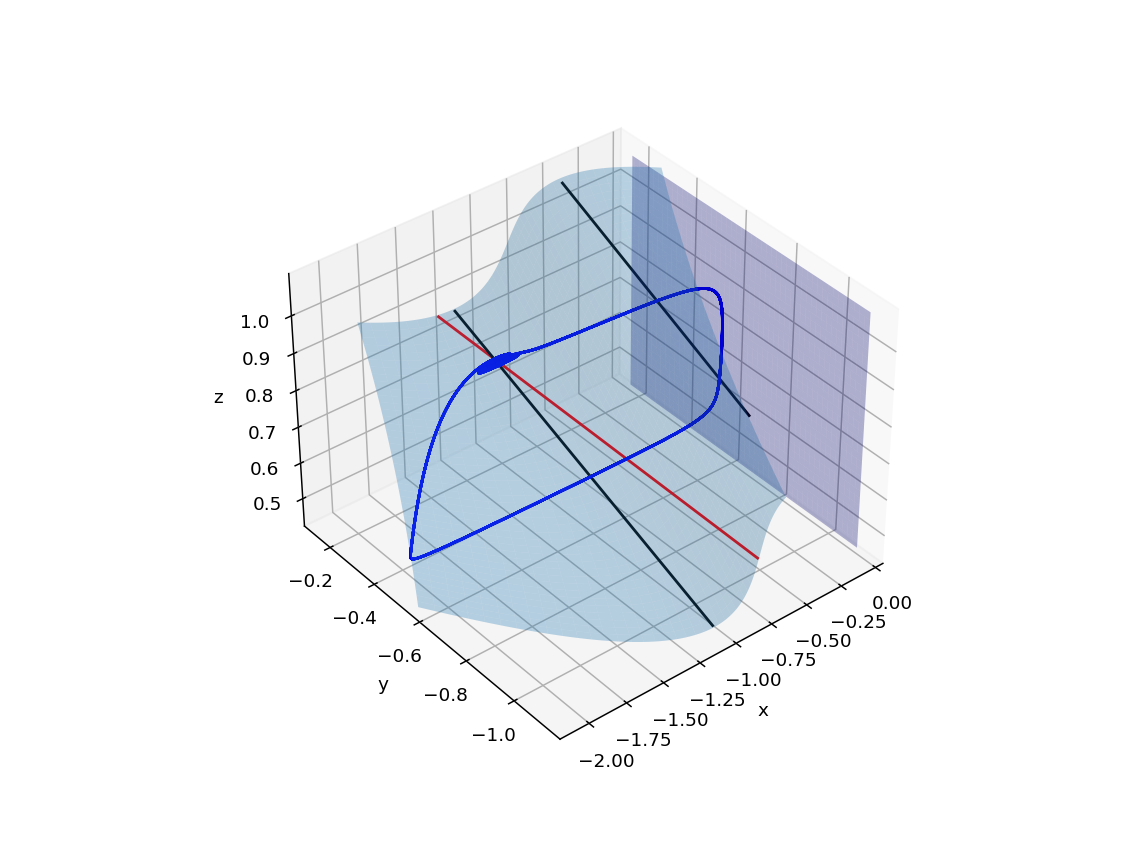}
	\caption{$c = 1.4$, $k=0.2$, $a=4.4$}
\end{subfigure}
	\caption{Given $ (c,k)\in \mc{V}_1 $, see \figref{regions}, oscillatory trajectories of Equation~\eqref{elno} feature plateaus above. For fixed $ (c,k)\in \mc{V}_1 $, there exists a unique $ a^- = a^-(c,k)>0 $ for which the equilibrium $\hat{p}_0$ of the reduced flow on $\mc{M}_2$ coincides with $ q^- $. For $ 0<a<a^-+\mc{O}(\delta,\rho) $, \eqref{elno} exhibits oscillatory dynamics. Moreover, for $ a $ close to $ a^- $, MMO trajectories with SAOs below are observed.}
	\figlab{V1}
\end{figure}

Finally, we observe that during the transition from oscillatory dynamics to steady state, Equation~\eqref{elno} features MMO trajectories with SAOs below; see \figref{V1}. Unfortunately, we are not able to take the approach described in \cite{kaklamanos2022bifurcations}, where we approximated the slow drift in order to predict the transition from relaxation oscillation to mixed-mode dynamics with SAOs below in dependence of $ a $: that drift cannot be deduced from the reduction in \eqref{reducedSsing} here due to $ z $ therein being the intermediate, rather than the slow, variable. 
{
\begin{remark}
While the slow drift could be approximated in the standard form of Equation~\eqref{reducedSz} after projection into the $(x,y)$-plane, cf.~\remref{standard}, preliminary analysis indicates that complications will arise due to the singular geometry of \eqref{elno}. A more in-depth investigation is included in plans for future work.
\end{remark}}

\subsubsection{\underline{$ (c,k) \in \mc{V}_2 = \mc{D}_2\cap\mc{A}_1$}}
Fix  $ (c,k) \in \mc{V}_2 $, as shown in \figref{regions}.

By \corref{Ds}, we have that $ P(q^-) \in \mc{S}^{a^+}$ and $ P(q^*) \in \mc{P}^a $. In addition, by \eqref{amp}, there exists $ a^- = a^-(c,k) $ such that $\hat{p}_0 \equiv q^- $, with $ \hat{p}_0\in \mc{S}^r $ for $ a<a^- $. Moreover, there exists $ a^+ = a^+(c,k) $ for which $\hat{p}_0 \equiv q^+ $, i.e. for which an equilibrium found by solving \eqref{eqbria} coincides with the folded singularity $ q^+ $. (However, we reiterate that for $ (c,k) $-values below the dashed purple curve $\mc{C}$ in \figref{regions}, the value $ a^+ $ is irrelevant, since trajectories cannot reach either ${q^-}$ or $\hat{p}_0$ in that case due to $x_{q^-}$ and $\hat{x}_0$ being positive.) Therefore, Equation~\eqref{elno} features both singular cycles that have intermediate segments in $\mc{M}_{\mc{S}}$ only and those that have segments in $\mc{M}_{\mc{P}}$ only.

Hence, given $ (c,k) \in\mc{V}_2$ and $ \delta,\rho>0$ sufficiently small, \eqref{elno} will exhibit oscillatory dynamics for $ a\in (a^+,a^-)+\mc{O}(\delta,\rho) $, respectively for $ 0<a<a^-+\mc{O}(\delta,\rho) $, if $ (c,k) $ lies above, respectively below, the dashed purple curve $\mc{C}$ in \figref{regions}, with trajectories projecting either onto the attracting portion $\mc{P}^a$ of the invariant plane $\mc{M}_{\mc P}$ or onto the attracting sheet $\mc{S}_{\delta\rho}^{a+}$ of the slow manifold under the fast flow of \eqref{elno}. Correspondingly, plateauless MMO trajectories will be observed in the latter scenario, whereas plateaus above will occur in the former. 

\begin{figure}[ht!]
	\centering
	\begin{subfigure}[b]{0.48\textwidth}
	\centering
	\begin{tikzpicture}
	\node at (0,0){
	\includegraphics[scale=0.25]{pics/elno/c14k4a3.png}};
	
	\node at (0.15,-2.9) {$t$};
	\node at (-2.9,0) {$x$};
	
	\end{tikzpicture}
		\caption{$c = 1.4$, $k=0.4$, $a=3$}
	\end{subfigure}
	\begin{subfigure}[b]{0.48\textwidth}
	\centering
	\begin{tikzpicture}
	\node at (0,0){
	\includegraphics[scale=0.25]{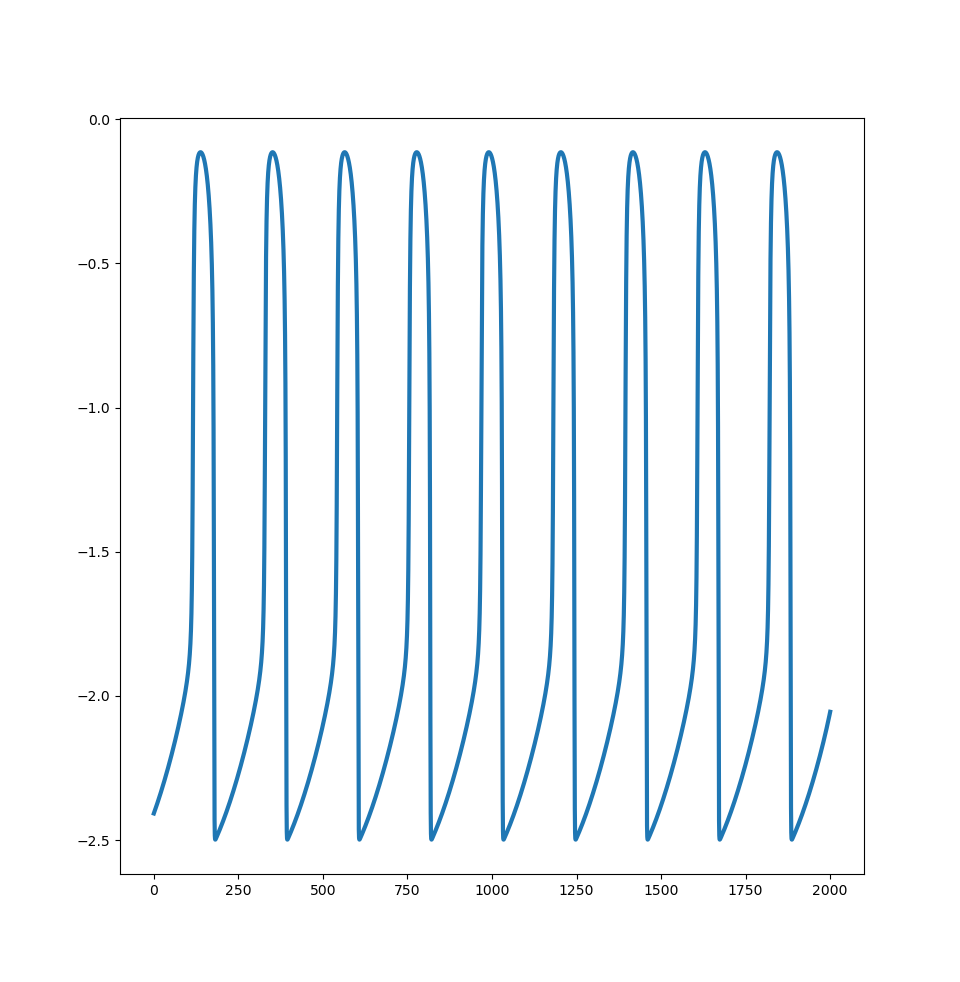}};
	
	\node at (0.15,-2.9) {$t$};
	\node at (-2.9,0) {$x$};
	
	\end{tikzpicture}
		\caption{$c = 1.4$, $k=0.4$, $a=20$}
	\end{subfigure}
	\\
	\begin{subfigure}[b]{0.48\textwidth}
		\centering
		\includegraphics[scale=0.2]{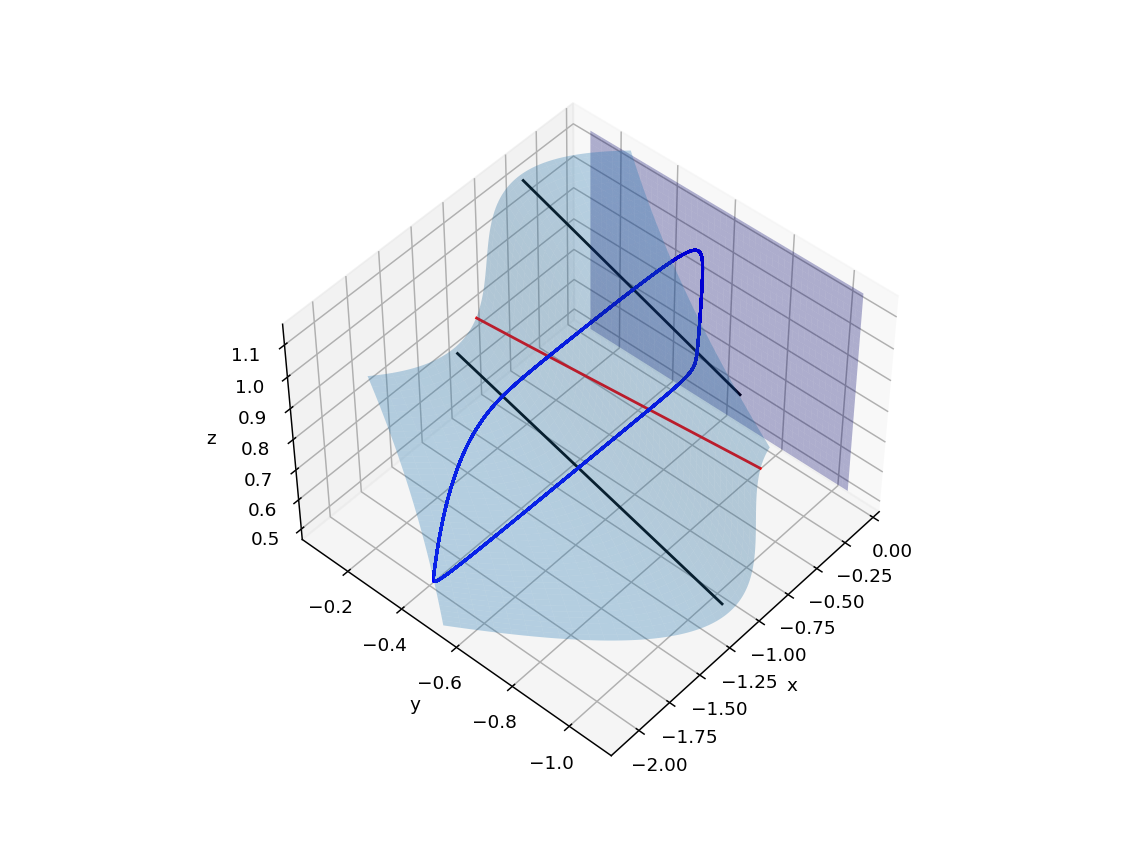}
		\caption{$c = 1.4$, $k=0.4$, $a=3$}
	\end{subfigure}
	\begin{subfigure}[b]{0.48\textwidth}
		\centering
		\includegraphics[scale=0.2]{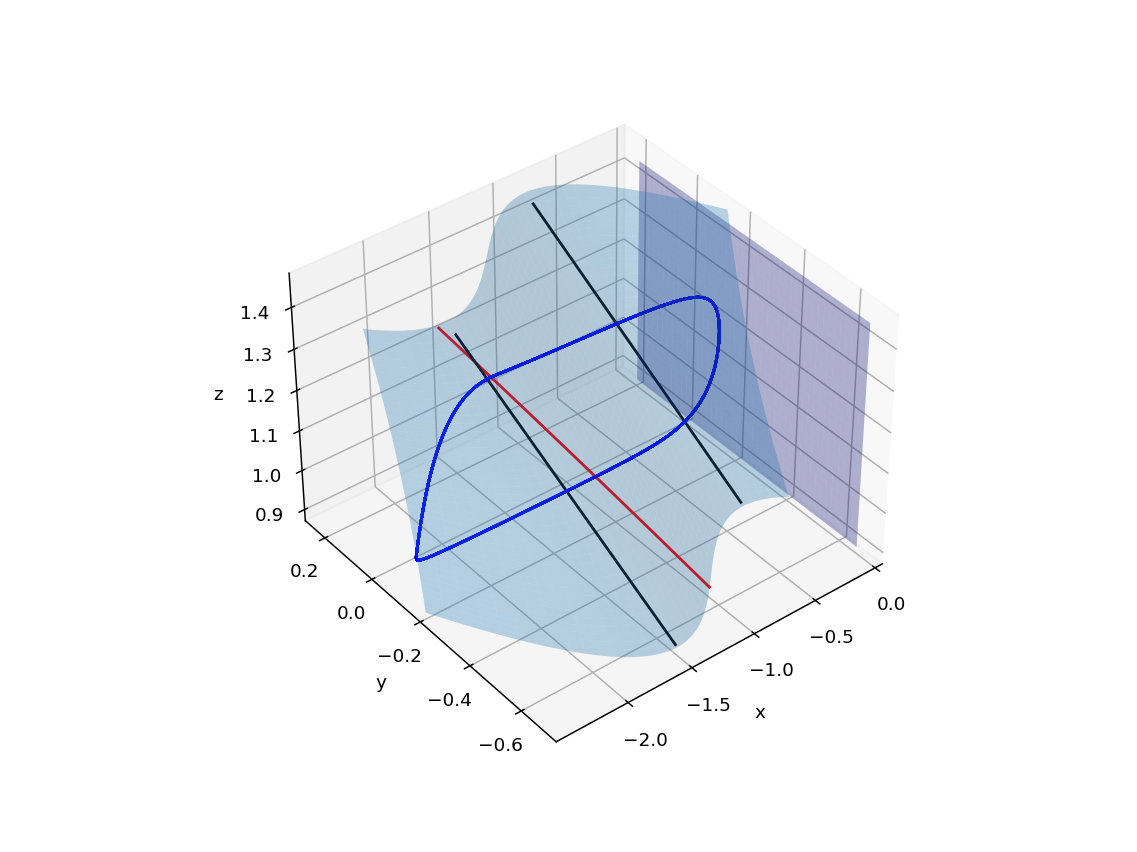}
		\caption{$c = 1.4$, $k=0.4$, $a=20$}
	\end{subfigure}
	
	\caption{Given $ (c,k)\in \mc{V}_2 $ below the dashed purple curve $\mc{C}$ in \figref{regions}, oscillatory trajectories of Equation~\eqref{elno} can either feature plateaus above or be plateauless. For fixed $ (c,k)\in \mc{V}_2 $, there exists a unique $ a^- = a^-(c,k)>0 $ for which an equilibrium of \eqref{elno} coincides with $ q^- $; \eqref{elno} exhibits oscillatory dynamics when $ 0<a<a^-+\mc{O}(\delta,\rho) $. Moreover, there exists $ a_p $ such that for $ a\lesssim a_p $, \eqref{elno} exhibits oscillation with plateaus above, whereas for $ a\gtrsim a_p $, plateauless trajectories are observed.}
	\figlab{V2}
\end{figure}

The above two scenarios are realised in dependence of the location of the global equilibrium $\hat{p}$ of \eqref{elno} for $\delta,\rho>0$ sufficiently small; the transition between the two hence depends on the value of the parameter $a$: as described in \secref{singgeomElno}, there exists an $a$-value $ a_p\in (a^+,a^-) $ such for $ a<a_p $, $ P(\mc{L}^-\cap \mc{W}(\hat{p}_0)) \in \mc{P}^a$, whereas $ P(\mc{L}^-\cap \mc{W}(\hat{p}_0)) \in \mc{S}^{a^+}$ for $ a>a_p $. However, we reiterate that the corresponding transition between trajectories with plateaus and those without seems continual rather than abrupt; a qualitative explanation is as follows: as $a$ is increasing and approaching the value $a_p=a_p(c,k)$, the flow of \eqref{elno} enters the vicinity of $\mc{P}^a$ closer and closer to $\mc{F}_\mc{P}$ in the $z$-direction, which implies that the segment along which the corresponding trajectory undergoes delayed loss of stability after crossing $\mc{F}_\mc{P}$ decreases until $a$ is such that the trajectory is attracted to $\mc{S}^{a^+}_{\delta\rho}$; cf. \figref{V2}.

\begin{remark}
    We remark that the distinction between relaxation oscillation with and without plateaus above is also relevant in the two-timescale context of \eqref{elno} with $\delta>0$ sufficiently small and $\rho=\mc{O}(1)$.
\end{remark}
Finally, fixing again $ (c,k) = (1.4,0.4) $ -- which lies below the dashed purple curve $\mc{C}$ in \figref{regions} -- we calculate that $ a^- \simeq 25.55  $. Numerically, we obtain that the flow of \eqref{elno} converges to steady state for $ a\gtrsim 25.5 $. Simulated trajectories for $ (c,k) $ fixed and various values of $ a $, with $ \delta=0.01=\rho $, are illustrated in \figref{V2}.

\subsubsection{\underline{$ (c,k) \in \mc{V}_3 = \mc{D}_3\cap\mc{A}_1$}}
Fix  $ (c,k) \in \mc{V}_3 $, as shown in \figref{regions}.

By \corref{Ds}, it holds that $ P(q^-), P(q^*) \in \mc{S}^{a^+}$. In addition, by \eqref{amp}, there exist $ a^- = a^-(c,k)>0 $ and $ a^+ = a^+(c,k)>0 $ for which $\hat{p}_0 \equiv q^- $ and $\hat{p}_0 \equiv q^+ $, respectively. It then follows that, for $ a\in(a^+,a^-) $, we have $ \hat{p}_0\in \mc{S}^r $. Hence, the location of the singularity $q^+$ is such that there can exist singular cycles which pass through $q^+$ and which feature segments on $\mc{S}^{a^\mp}$. Since, moreover, the reduced flow on $\mc{M}_{2\mc{S}}\cap\mc{S}^{a^\mp}$ is directed towards $q^\mp$, respectively, singular cycles of \eqref{elno} contain no segments on $\mc{M}_{\mc{P}}$ and are hence plateauless. 

Therefore, for $ \delta,\rho>0$ sufficiently small, Equation~\eqref{elno} features oscillatory dynamics for $a\in (a^+,a^-)+\mc{O}(\delta,\rho) $, where we emphasise that trajectories with plateaus are not possible.  

In particular, fixing $ (c,k) = (1.06,0.4) $, we calculate that $ a^+ \simeq 0.2$ and $ a^- \simeq 61.02  $. Numerically, we obtain that the flow of \eqref{elno} converges to steady state for $ a\lesssim 0.5 $, as well as for $ a\gtrsim 40 $. (This seeming discrepancy for large $ a $-values is addressed in \remref{breakdown}.) Simulated trajectories for $ (c,k) $ fixed in this regime and various values of $ a $, with $ \delta=0.01=\rho $, are illustrated in \figref{V3}.

\begin{figure}[ht!]
	\centering
	\begin{subfigure}[b]{0.48\textwidth}
	\centering
	\begin{tikzpicture}
	\node at (0,0){
	\includegraphics[scale=0.25]{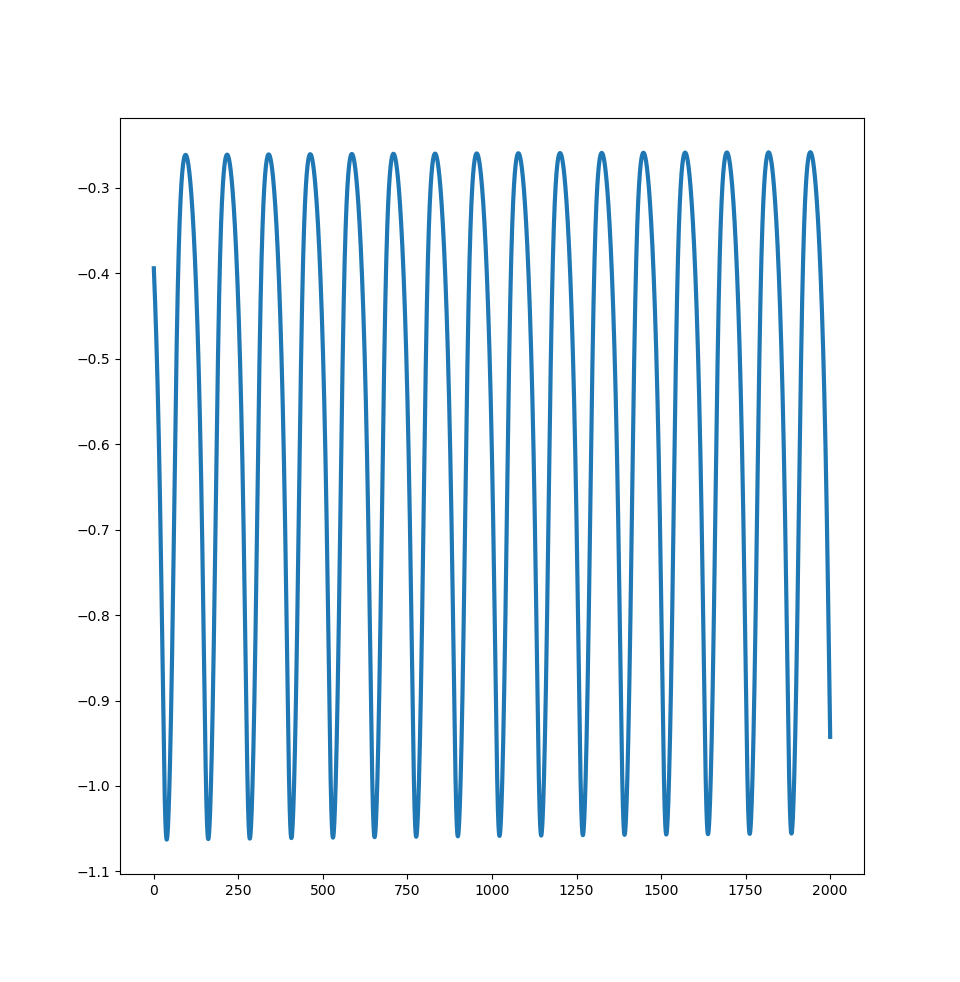}};
	
	\node at (0.15,-2.9) {$t$};
	\node at (-2.9,0) {$x$};
	
	\end{tikzpicture}
	\caption{$c = 1.06$, $k=0.4$, $a=1$}
	\end{subfigure}
	\begin{subfigure}[b]{0.48\textwidth}
	\centering
	\begin{tikzpicture}
	\node at (0,0){
	\includegraphics[scale=0.25]{pics/elno/c106k4a5.png}};
	
	\node at (0.15,-2.9) {$t$};
	\node at (-2.9,0) {$x$};
	
	\end{tikzpicture}
	\caption{$c = 1.06$, $k=0.4$, $a=5$}
	\end{subfigure}
	\\
	\begin{subfigure}[b]{0.48\textwidth}
		\centering
		\includegraphics[scale=0.2]{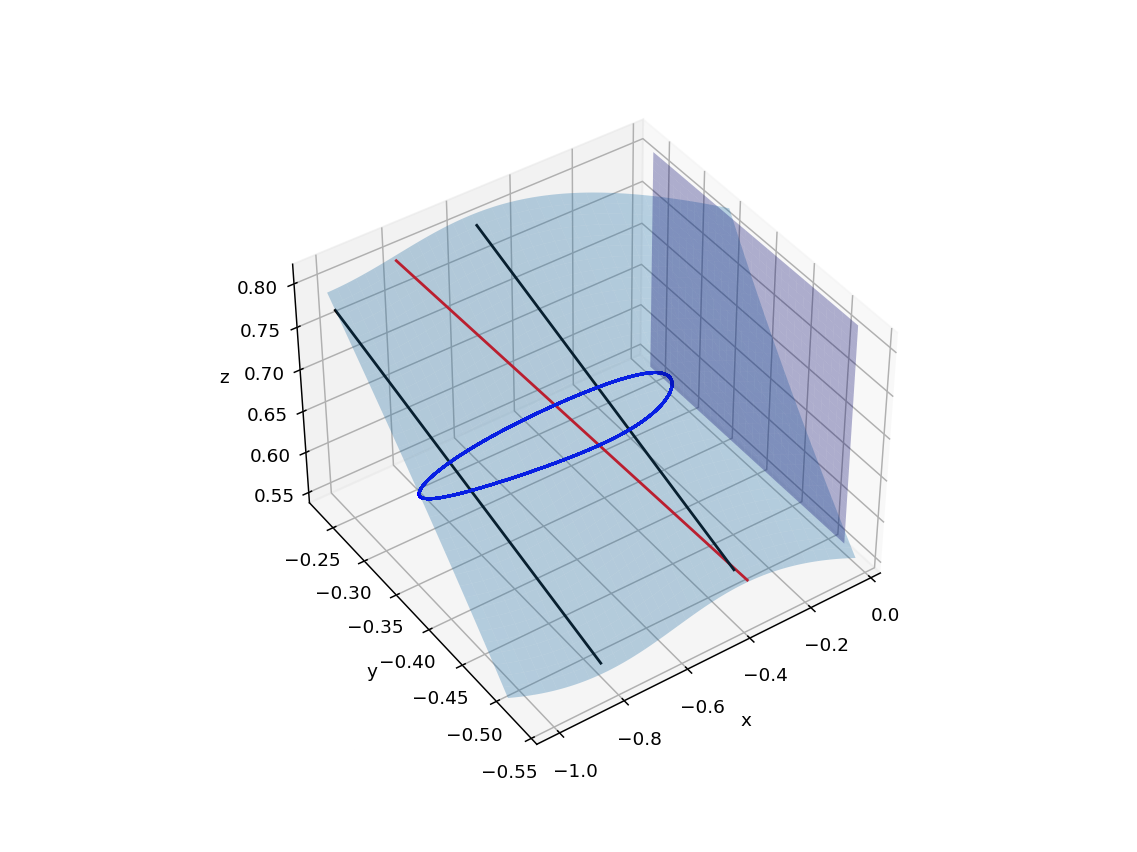}
		\caption{$c = 1.06$, $k=0.4$, $a=1$}
	\end{subfigure}
	\begin{subfigure}[b]{0.48\textwidth}
		\centering
		\includegraphics[scale=0.2]{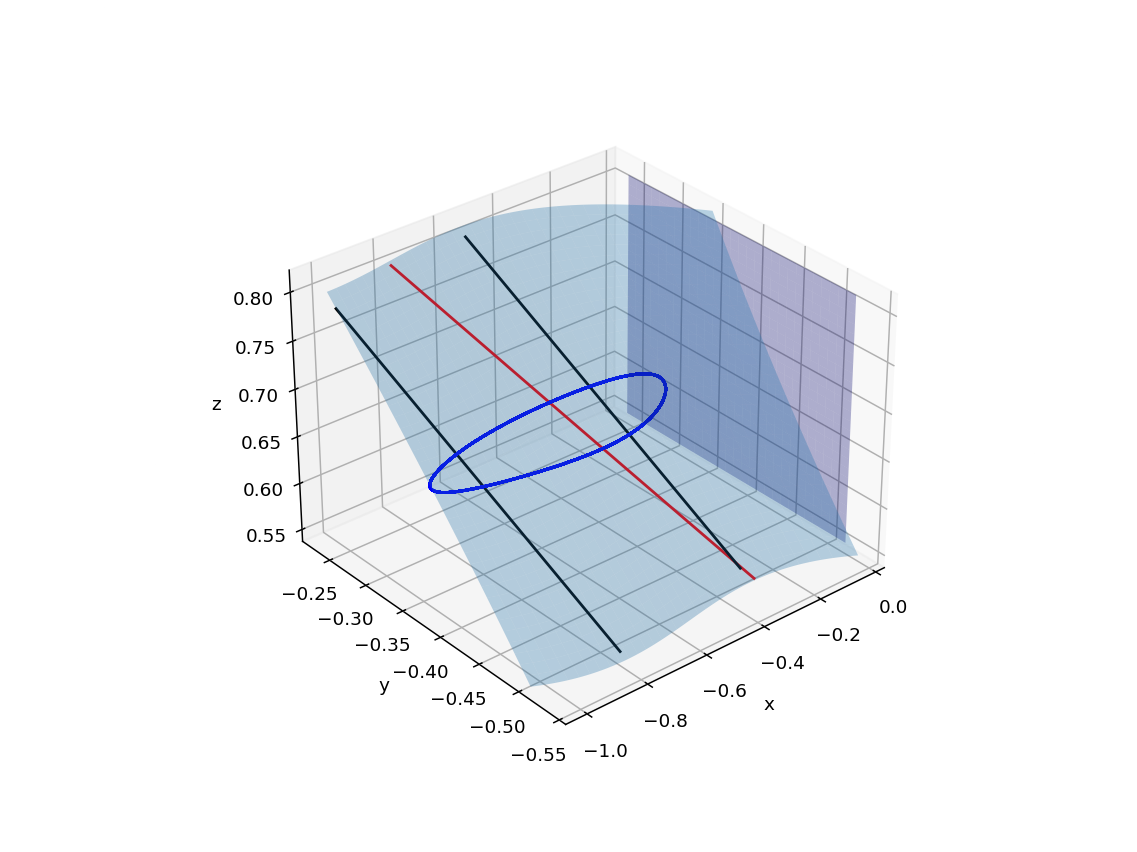}
		\caption{$c = 1.06$, $k=0.4$, $a=5$}
	\end{subfigure}
	
	\caption{Given $ (c,k)\in \mc{V}_3 $, see \figref{regions}, oscillatory trajectories of Equation~\eqref{elno} cannot feature plateaus above. For fixed $ (c,k)\in \mc{V}_3 $, there exist unique $a$-values $ a^- = a^-(c,k)>0 $ and $ a^+ = a^+(c,k)>0 $ for which an equilibrium of \eqref{elno} coincides with $ q^- $ and $ q^+ $, respectively. For $ a\in (a^+,a^-)+\mc{O}(\delta,\rho) $, \eqref{elno} exhibits oscillatory dynamics.}
	\figlab{V3}
\end{figure}

\subsubsection{\underline{$ (c,k) \in \mc{V}_4 = \mc{D}_2\cap\mc{A}_2$}}
Fix  $ (c,k) \in \mc{V}_4 $, as shown in \figref{regions}.

By \corref{Ds}, we have $ P(q^-) \in \mc{S}^{a^+}$ and $ P(q^*) \in \mc{P}^a $. In addition, by \eqref{amp}, there exists $ a^+ = a^+(c,k)>0 $ for which $\hat{p}_0 \equiv q^+ $, with $ \hat{p}_0\in \mc{S}^r $ for $ a>a^+ $; that value of $a$ again to leading order indicates a (singular) Hopf bifurcation of the perturbed system in \eqref{elno}, with $\delta,\rho>0$ sufficiently small -- we emphasise that there exists \textit{no} $ a^- = a^-(c,k)>0 $ for which $\hat{p}_0 \equiv q^- $. Hence, if $ (c,k) $ lies above the dashed purple curve $\mc{C}$ in \figref{regions}, then \eqref{elno} is expected to feature oscillatory dynamics for $a>a^++\mc{O}(\delta,\rho)$, whereas if $ (c,k) $ lies below that curve, then \eqref{elno} is expected to feature oscillatory dynamics for $a>0$.

In particular, fixing $ (c,k) = (1.4,0.7) $ -- which lies above the dashed purple curve $\mc{C}$ in \figref{regions} -- we calculate that $ a^+ \simeq 0.1$. Numerically, we obtain that the flow of \eqref{elno} converges to steady state for $ a\lesssim 0.1 $. Simulated sample trajectories for $ (c,k) $ fixed in this regime and various values of $ a $, with $ \delta=0.01=\rho $, are illustrated in \figref{V4}.

We remark that, similarly to the regime where $ (c,k) \in \mc{V}_3 $, the above analysis is valid only for $ a = \mc{O}(1) $: when $ a\gg1 $, the three-timescale separation where $ x $ is fast, $ y $ is intermediate, and $ z $ is slow, as described here, breaks down in accordance with \remref{breakdown}. An in-depth study of that case is again included in plans for future work. 

Moreover, as in the regime where $ (c,k) \in \mc{V}_2 $, we observe that there exists a value $ a_p = a_p(c,k) $ which distinguishes between trajectories that feature plateaus from those that do not; specifically, \eqref{elno} features oscillation with plateaus above for $ 0<a\lesssim a_p $, 
whereas plateauless trajectories are observed for $ a\gtrsim a_p $ and $ a = \mc{O}(1) $.

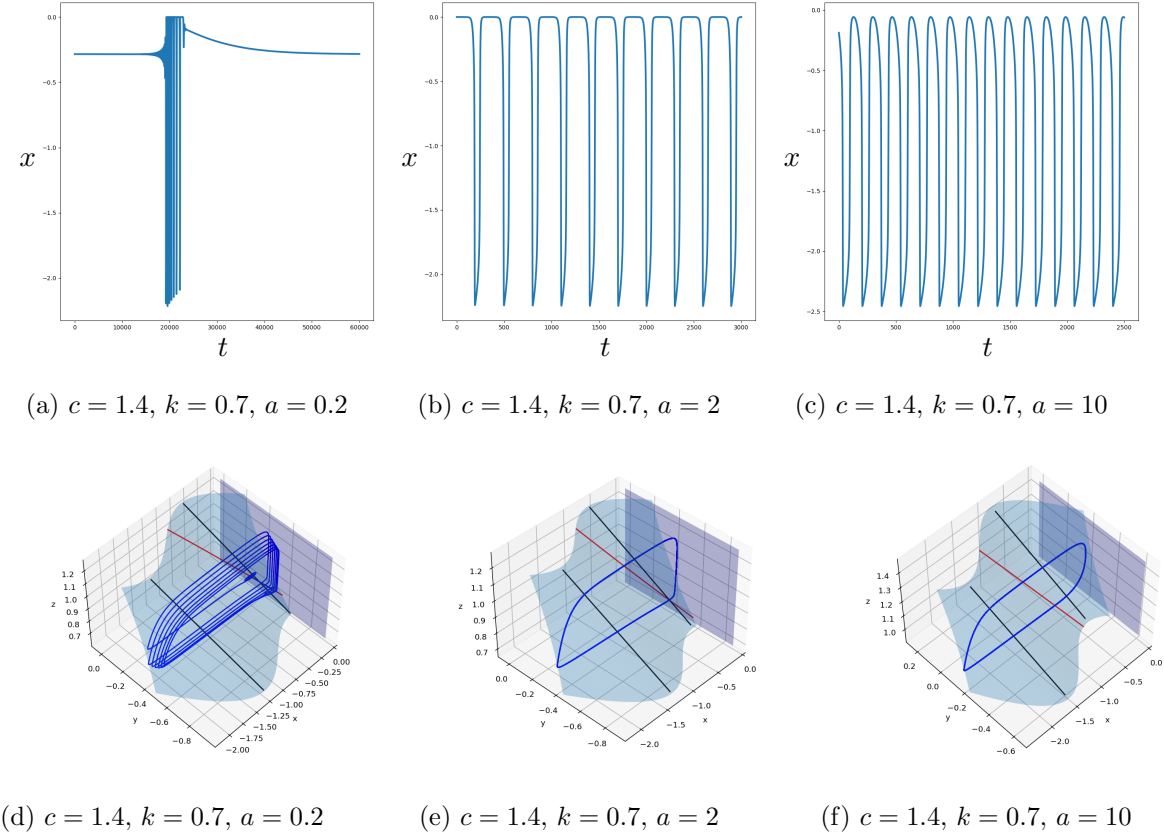
\begin{figure}[ht!]
	\centering
	\begin{subfigure}[b]{0.3\textwidth}
	\centering
	\begin{tikzpicture}
	\node at (0,0){
	\includegraphics[scale=0.22]{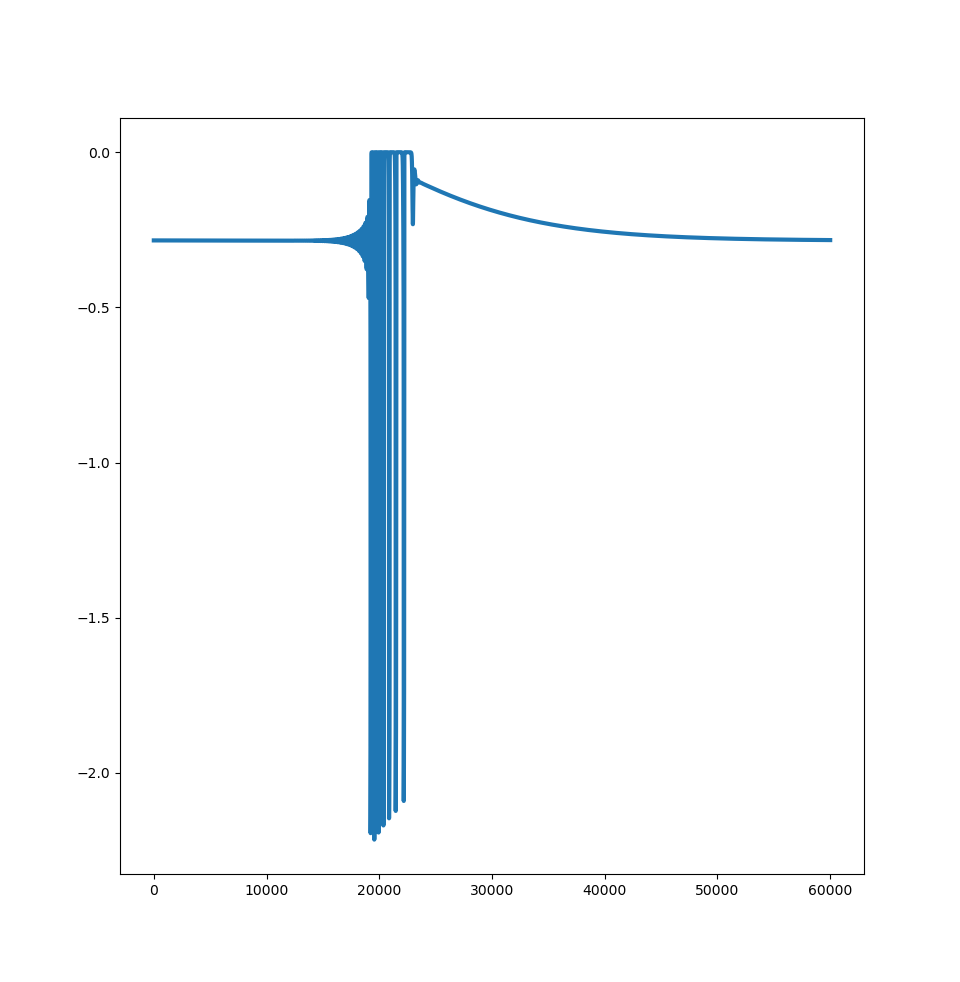}};
	
	\node at (0.15,-2.5) {$t$};
	\node at (-2.45,0) {$x$};
	
	\end{tikzpicture}
	\caption{$c = 1.4$, $k=0.7$, $a=0.2$}
	\end{subfigure}
	\begin{subfigure}[b]{0.3\textwidth}
	\centering
	\begin{tikzpicture}
	\node at (0,0){
	\includegraphics[scale=0.22]{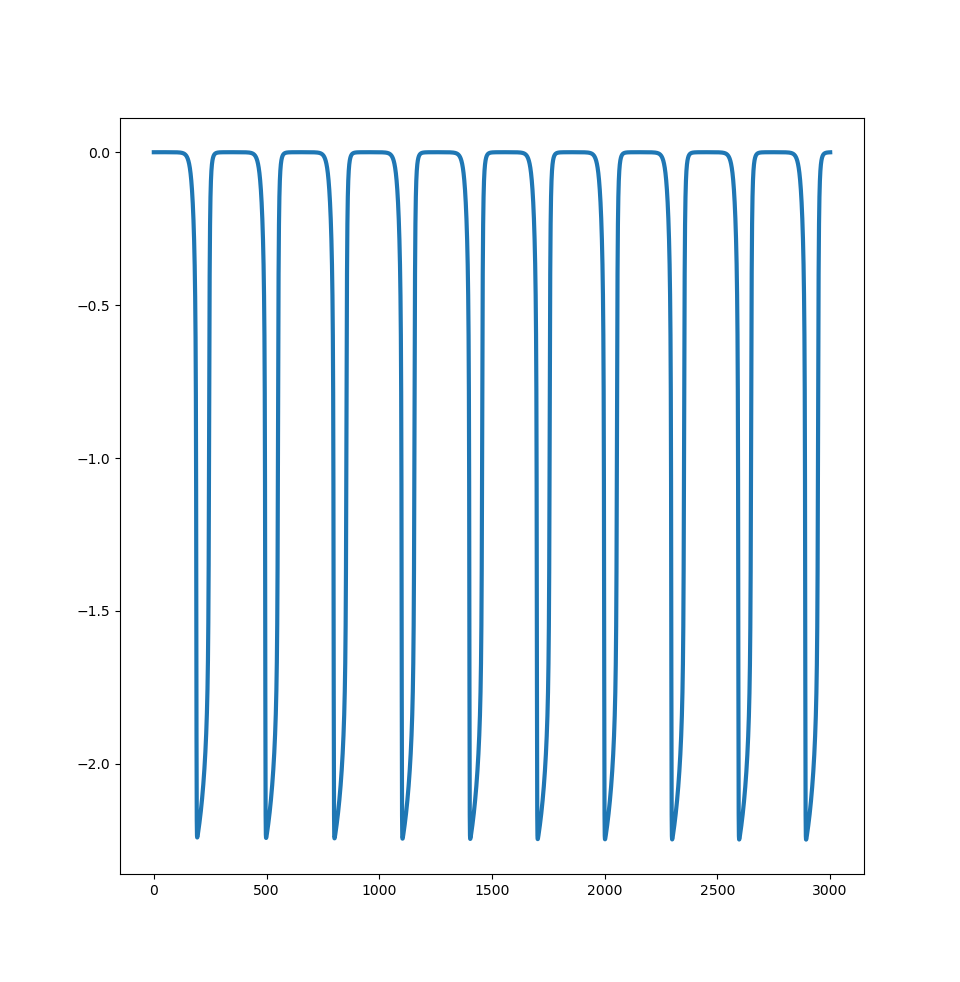}};
	
	\node at (0.15,-2.5) {$t$};
	\node at (-2.45,0) {$x$};
	
	\end{tikzpicture}
	\caption{$c = 1.4$, $k=0.7$, $a=2$}
	\end{subfigure}
	\begin{subfigure}[b]{0.3\textwidth}
	\centering
	\begin{tikzpicture}
	\node at (0,0){
	\includegraphics[scale=0.22]{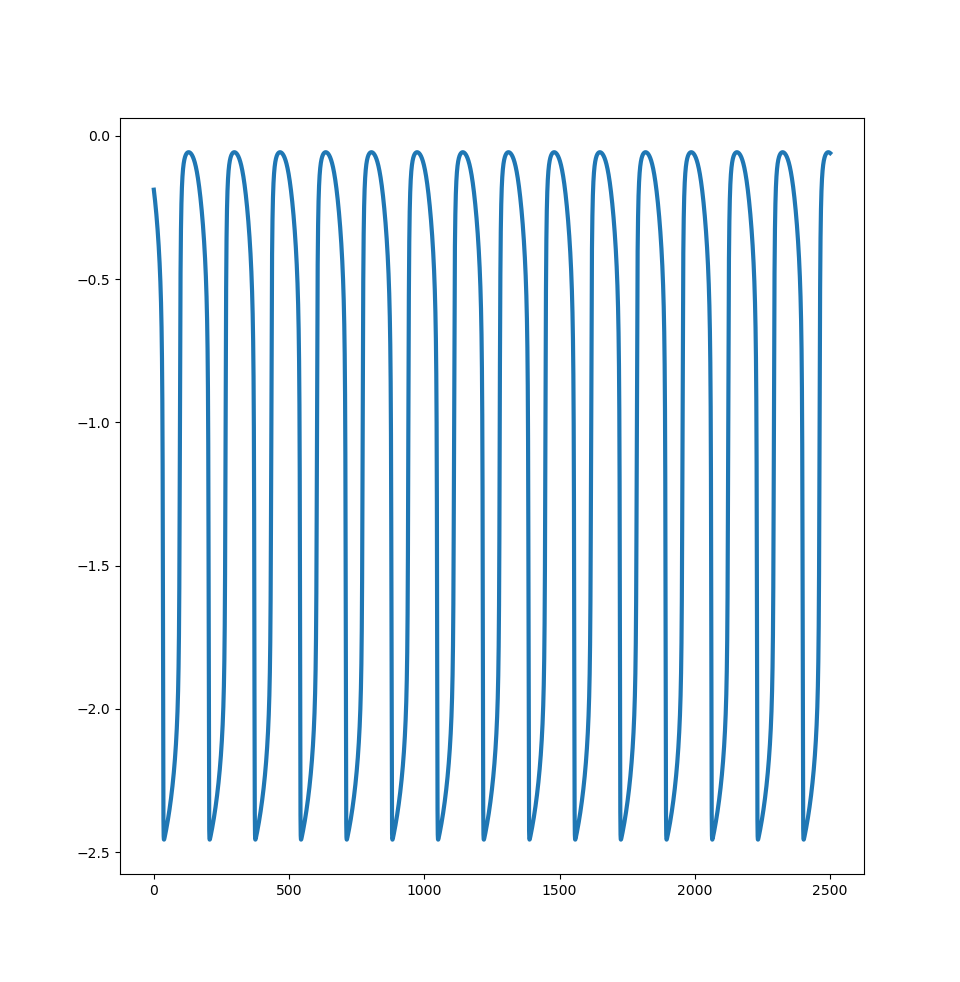}};
	
	\node at (0.15,-2.5) {$t$};
	\node at (-2.45,0) {$x$};
	
	\end{tikzpicture}
	\caption{$c = 1.4$, $k=0.7$, $a=10$}
	\end{subfigure}
	\\
	\begin{subfigure}[b]{0.32\textwidth}
		\centering
		\includegraphics[scale=0.16]{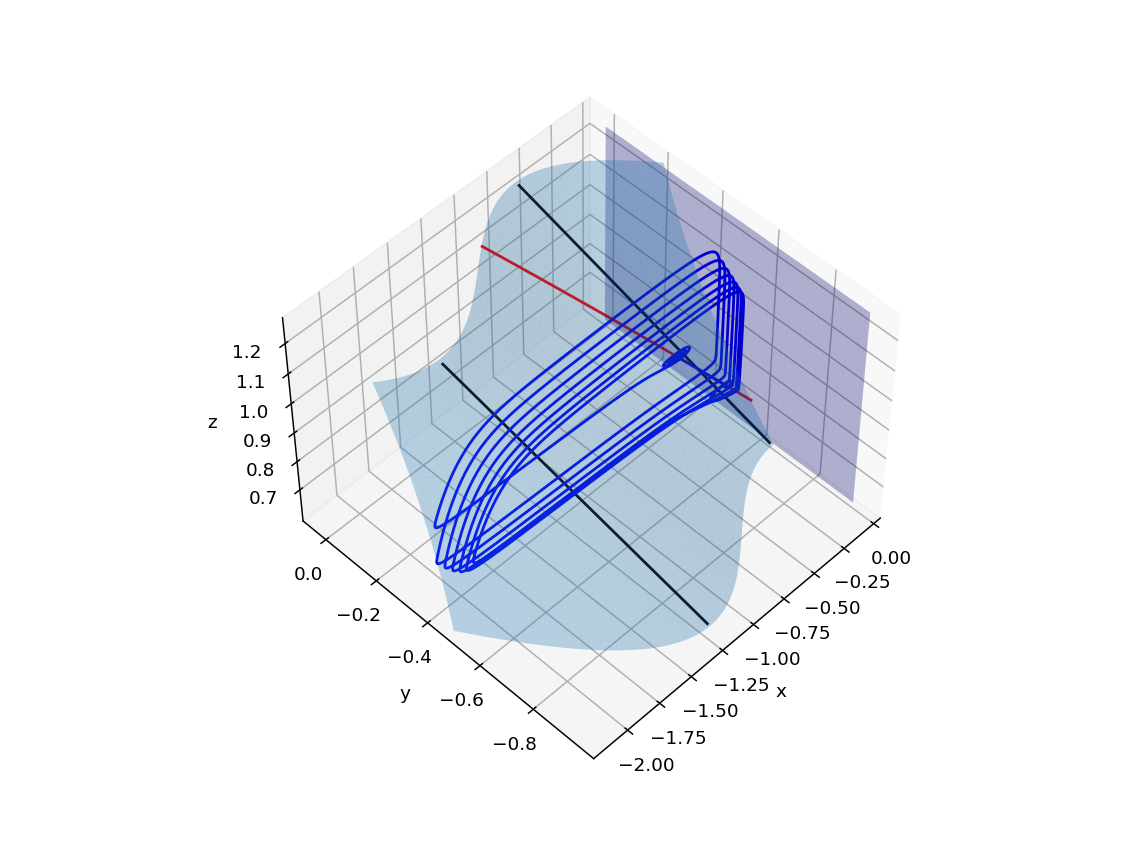}
		\caption{$c = 1.4$, $k=0.7$, $a=0.2$}
	\end{subfigure}
	\begin{subfigure}[b]{0.32\textwidth}
		\centering
		\includegraphics[scale=0.16]{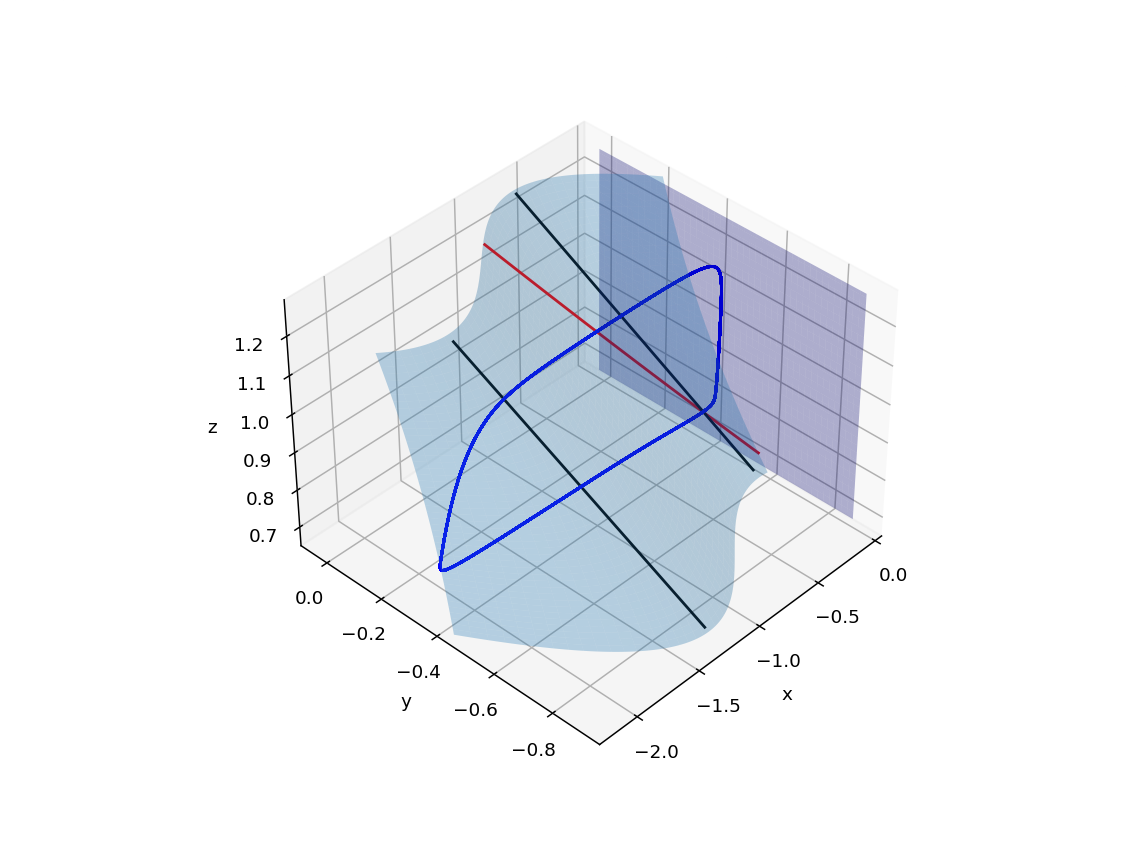}
		\caption{$c = 1.4$, $k=0.7$, $a=2$}
	\end{subfigure}
	\begin{subfigure}[b]{0.32\textwidth}
		\centering
		\includegraphics[scale=0.16]{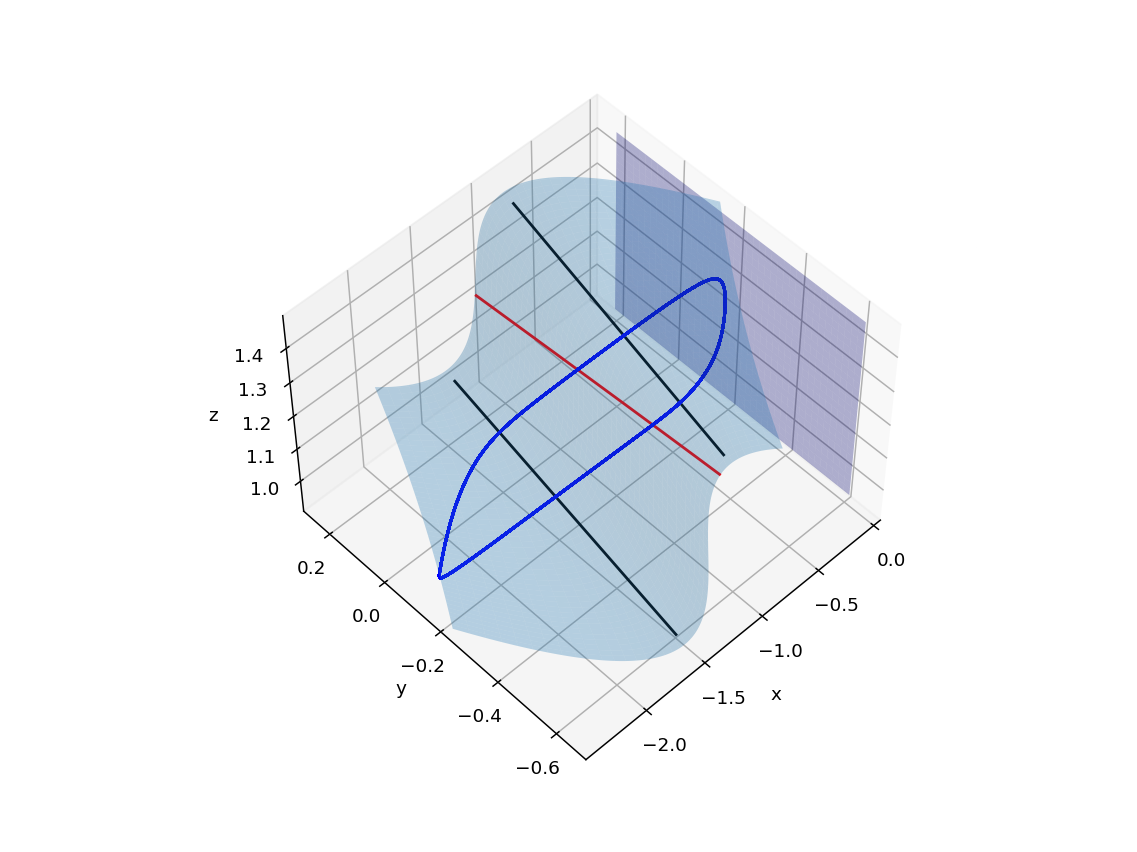}
		\caption{$c = 1.4$, $k=0.7$, $a=10$}
	\end{subfigure}
	\caption{Given $ (c,k)\in \mc{V}_4 $, see \figref{regions}, oscillatory trajectories of Equation~\eqref{elno} can either feature plateaus above or be plateauless. For fixed $ (c,k)\in \mc{V}_4 $, there exists a unique $ a^+ = a^+(c,k)>0 $ for which an equilibrium point of \eqref{elno} coincides with  $ q^+ $. For $ a>a^++\mc{O}(\delta,\rho) $, Equation~\eqref{elno} exhibits oscillatory dynamics. Moreover, there exists $ a_p $ such that for $ a\lesssim a_p $, \eqref{elno} features oscillatory trajectories with plateaus above, while for $ a\gtrsim a_p $, plateauless oscillation is observed. Finally, for $ a $-values close to $a^+ $, MMOs with SAOs above occur; note that the LAO component of the latter contains plateaus due to the corresponding intermediate segments evolving on $ \mc{P}^{a} $; cf.~\figref{V5} for comparison. }
	\figlab{V4}
\end{figure}

\subsubsection{\underline{$ (c,k) \in \mc{V}_5 = \mc{D}_3\cap\mc{A}_2$}}
Fix  $ (c,k) \in \mc{V}_5 $, as shown in \figref{regions}.

By \corref{Ds}, it holds that $ P(q^-), P(q^*) \in \mc{S}^{a^+}$. In addition, by \eqref{amp}, there exists $ a^+ = a^+(c,k)>0 $ for which $\hat{p}_0 \equiv q^+ $; for $ a>a^+ $, it holds that $ \hat{p}_0\in \mc{S}^r $ -- we emphasize that there again exists \textit{no} $ a^- = a^-(c,k)>0 $ for which $\hat{p}_0 \equiv q^- $. Therefore, for $ \delta,\rho>0$ sufficiently small, the system in \eqref{elno} features oscillatory dynamics for $ a> a^++\mc{O}(\delta,\rho) $ and $ a = \mc{O}(1) $, as in the regime where $(c,k)\in\mc{V}_4$. However, the corresponding MMO trajectories are again plateauless, which can be reasoned as in the regime where $(c,k)\in\mc{V}_3$.

In particular, fixing $ (c,k) = (1.2,0.7) $, we calculate that $ a^+ \simeq 1.6$. Numerically, we obtain that the flow of \eqref{elno} converges to steady state for $ a\lesssim 2.2 $. Simulated sample trajectories for $ (c,k) $ fixed in this regime and various values of $ a $, with $ \delta=0.01=\rho $, are illustrated in \figref{V5}.

\begin{figure}[ht!]
	\centering
	\begin{subfigure}[b]{0.48\textwidth}
	\centering
	\begin{tikzpicture}
	\node at (0,0){
	\includegraphics[scale=0.25]{pics/elno/c12k7a2-2.png}};
	
	\node at (0.15,-2.9) {$t$};
	\node at (-2.9,0) {$x$};
	
	\end{tikzpicture}
	\caption{$c = 1.2$, $k=0.7$, $a=2.2$}
	\end{subfigure}
	\begin{subfigure}[b]{0.48\textwidth}
	\centering
	\begin{tikzpicture}
	\node at (0,0){
	\includegraphics[scale=0.25]{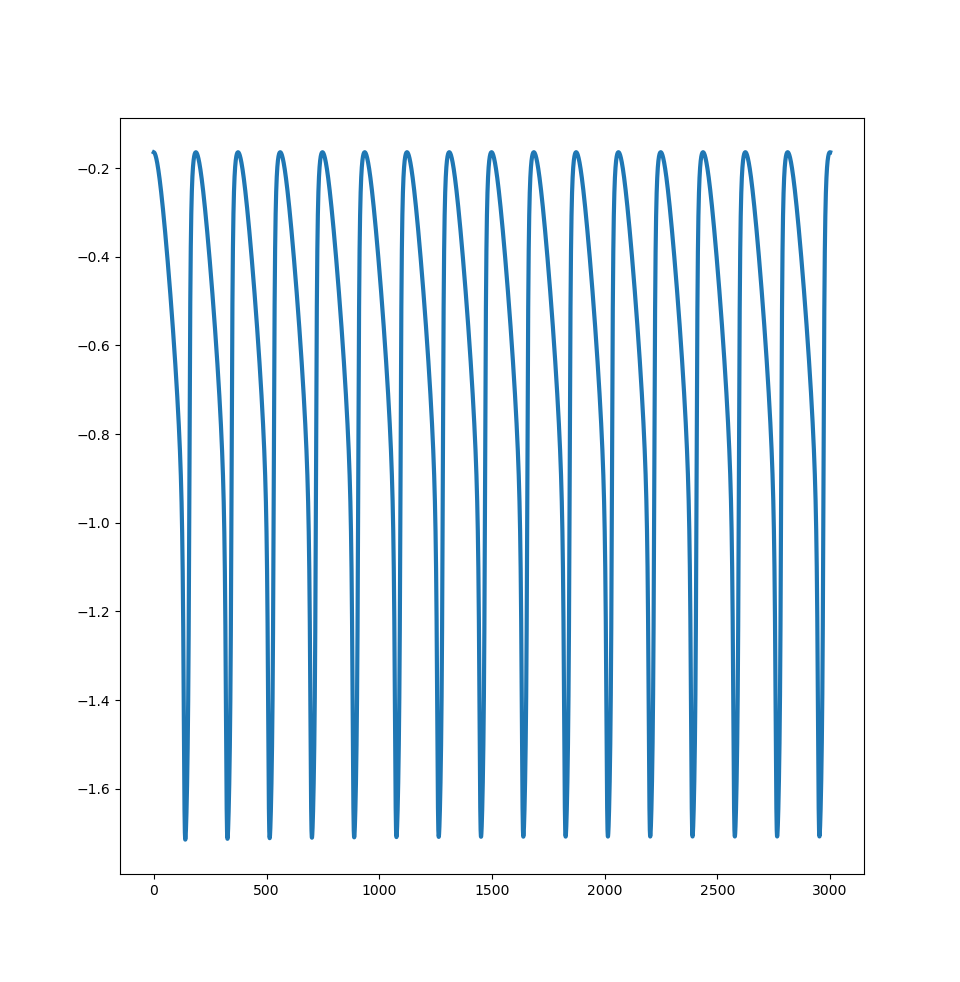}};
	
	\node at (0.15,-2.9) {$t$};
	\node at (-2.9,0) {$x$};
	
	\end{tikzpicture}
	\caption{$c = 1.2$, $k=0.7$, $a=3$}
	\end{subfigure}
	\\
	\begin{subfigure}[b]{0.48\textwidth}
		\centering
		\includegraphics[scale=0.2]{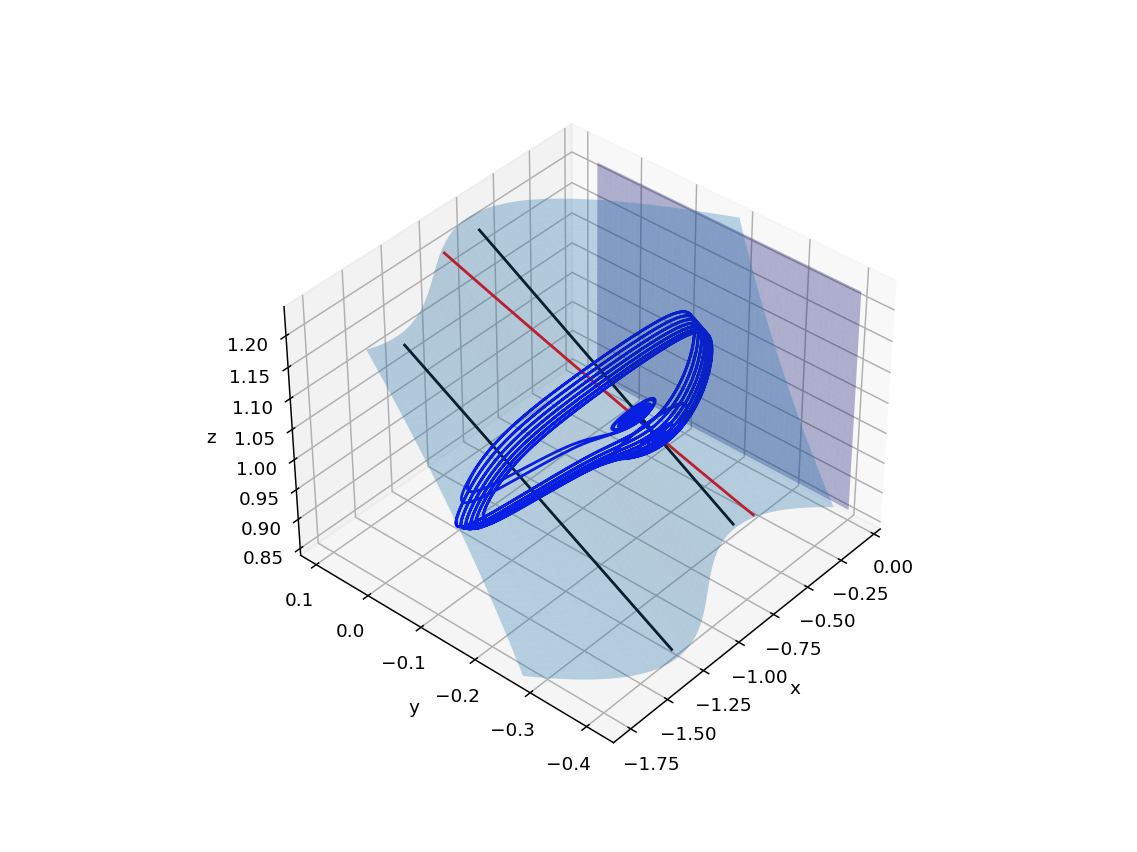}
		\caption{$c = 1.2$, $k=0.7$, $a=2.2$}
	\end{subfigure}
	\begin{subfigure}[b]{0.48\textwidth}
		\centering
		\includegraphics[scale=0.2]{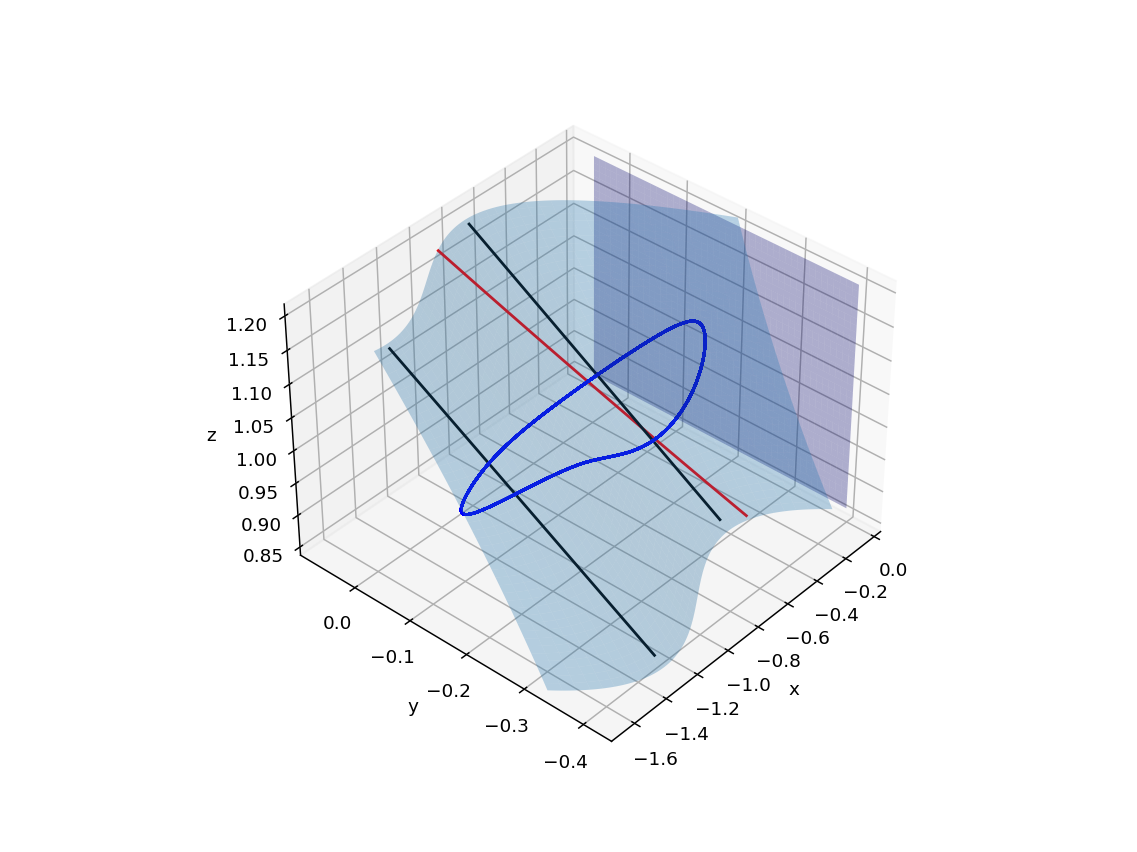}
		\caption{$c = 1.2$, $k=0.7$, $a=3$}
	\end{subfigure}
	
	\caption{Given $ (c,k)\in \mc{V}_5 $, see \figref{regions}, oscillatory trajectories of Equation~\eqref{elno} cannot feature plateaus above. For fixed $ (c,k)\in \mc{V}_5 $, there exists a unique $ a^+ = a^+(c,k)>0 $ for which an equilibrium point of \eqref{elno} coincides with $ q^+ $. For $ a>a^+ $, \eqref{elno} exhibits oscillatory dynamics. Finally, for $a$-values close to $a^+ $, MMOs with SAOs above occur; note that the LAO component of the latter consists of ``typical'' relaxation segments due to the corresponding intermediate segments evolving on $ \mc{S}^{a^-} $; cf.~\figref{V4} for comparison.}
	\figlab{V5}
\end{figure}

We remark that, similarly to the previous regimes, the above analysis is valid only for $ a = \mc{O}(1) $; for $ a\gg1 $, the $a$-dependent terms in \eqref{elno} become large, which implies that the three-timescale separation where $ x $ is fast, $ y $ is intermediate and $ z $ is slow as described here is no longer valid, in accordance with \remref{breakdown}.

\subsubsection{\underline{$ (c,k) \in \mc{V}_6 = \mc{D}_3\cap\mc{A}_3$}}
Fix $ (c,k) \in \mc{V}_6 $, as shown in \figref{regions}. Then, there exists no $ a>0 $ for which \eqref{elno} features oscillatory dynamics, since the relevant $ a $-values $a^\mp$ are defined only for $ (c,k) \in \mc{A}_1\cup \mc{A}_2 $; recall \corref{As}. Rather, the flow of \eqref{elno} converges to steady state in this regime. 

	

\section{Conclusion}
\seclab{elno-sum}

In this work, we have studied the dynamics of a three-dimensional ordinary differential equation model for the El Ni\~no Southern Oscillation (ENSO) phenomenon, by extending the analysis in \cite{roberts2016mixed} to the three-timescale context of Equation~\eqref{elno}, i.e. by considering $ \delta,\rho>0 $ sufficiently small therein. We have explored the properties of oscillatory trajectories in dependence of the parameters $ c,k $, and $ a $ in Equation~\eqref{elno}, and we have associated them with the geometric properties of the system in the singular limit of $ \delta=0=\rho $. 

To that end, in \secref{singgeomElno}, we studied the geometric properties of the critical and 2-critical manifolds for \eqref{elno} in that double singular limit: we showed that the resulting critical manifold $\mc{M}_1$ for \eqref{elno} is self-intersecting, in that it consists of an $S$-shaped portion $\mc{M}_{\mc S}$ which intersects transversely with a planar manifold $\mc{M}_{\mc P}$. Correspondingly, the $2$-critical manifold $\mc{M}_2$ consists of a portion $\mc{M}_{2\mc{S}}$ within $\mc{M}_{\mc S}$ which is again $S$-shaped, as well as of a linear portion $\mc{M}_{2\mc{P}}$ within $\mc{M}_{\mc{P}}$.

In particular, we illustrated various geometric configurations that result from the possible relative locations of these manifolds, as well as of sets thereon where normal hyperbolicity is lost, in dependence of the parameters $c$, $k$, and $a$ in \eqref{elno}. Specifically, we showed that the parameter $ c $ is associated with the geometric properties of the two-dimensional critical manifold $ \mc{M}_1 $ of \eqref{elno} in the singular limit of $ \delta=0 $. For fixed $ c>1 $, the parameter $ k\in(0,1) $ then determines the geometric properties of the 2-critical manifold $ \mc{M}_2 $ in the double singular limit of $ \delta=0=\rho $. Finally, the parameter $ a $ does not affect the singular geometry of Equation~\eqref{elno}; rather, given a fixed geometry, it can distinguish between steady-state behaviour and oscillatory dynamics, as well as between qualitatively different oscillatory behaviours, via the reduced flow on the corresponding invariant manifolds in \eqref{elno}. Crucially, we related our discussion of the geometry of \eqref{elno} for $\delta=0=\rho$ to the properties of the associated singular cycles.

In \secref{mains}, we illustrated various types of oscillatory dynamics in \eqref{elno} in dependence of the possible geometric configurations for $ \delta,\rho>0 $ sufficiently small, i.e. by perturbing off the singular picture constructed in \secref{singgeomElno}. By classifying the parameter regimes corresponding to the various geometric configurations that are observed, we were thus able to uncover novel dynamics that had not been documented in previous works. In particular, by reference to \figref{regions}, we have shown the following.
\begin{enumerate}
	\item If $ (c,k)\in\mc{V}_1 $, only oscillation with plateaus above is possible; cf.~\figref{V1}. 
	\item If $ (c,k)\in\mc{V}_2 $, oscillation with or without plateaus is possible, in dependence of the parameter $ a $; cf.~\figref{V2}.
	\item If $ (c,k)\in\mc{V}_3 $, only plateauless oscillation is possible; cf.~\figref{V3}.
	\item If $ (c,k)\in\mc{V}_4 $, oscillation with or without plateaus is possible, in dependence of the parameter $ a $; moreover, MMO trajectories with plateaus above can also potentially feature segments of SAOs above; cf.~\figref{V4}.
	\item If $ (c,k)\in\mc{V}_5 $, only plateauless oscillation is possible, which, in dependence of the parameter $ a $, can also potentially feature SAOs above; cf.~\figref{V5}.
	\item Finally, if $ (c,k)\in\mc{V}_6 $, no oscillatory dynamics is possible, with the flow of \eqref{elno} converging to steady state.
\end{enumerate}

{While we have hence painted a relatively complete picture of the qualitative dynamics of Equation~\eqref{elno}, a number of questions remain open for future investigation. First, and as alluded to repeatedly, it would be of interest to consider the scenario where $a\gg 1$ in \eqref{elno}, in which the separation of scales assumed throughout here breaks down. Second, the standard form introduced briefly in \remref{standard} could be explored further, in particular with regard to its utility for describing the transition from relaxation oscillation to mixed-mode dynamics with epochs of SAOs; on a related note, a precise characterisation of the latter would seem relevant. Finally, the emergence of canard trajectories along the repelling branch $\mc{S}_{\delta\rho}^r$ of $\mc{S}$ that can potentially be associated with the curvature of that manifold \cite{desroches2011canards} could be investigated; recall \figref{V5}.}

\section*{Acknowledgements}
The content of this work was part of the first author's PhD thesis, completed between 2018 and 2021 at the University of Edinburgh \cite{kaklamanos2022mixed}, and was motivated by a dissertation that had been submitted by Yichen Su for the degree of an MSc in Computational Applied Mathematics in 2020. The authors would like to thank Tom Mackay and Martin Wechselberger for feedback and recommendations that led to an improved version of the manuscript.

\bibliographystyle{siam}
\bibliography{refs}
\end{document}